# Some applications of the sine and cosine integrals

Donal F. Connon

dconnon@btopenworld.com

1 December 2012

**Abstract**

We show how the sine and cosine integrals may be usefully employed in the evaluation of some more complex integrals.

## 1. Introduction

As noted by Glaisher [23], the sine and cosine integrals were introduced by Schlömilch [35] in 1846 to evaluate integrals of the form

$$\int_0^\infty \frac{\sin(tx)}{a^2 - x^2} dx \text{ and } \int_0^\infty \frac{x\cos(tx)}{a^2 - x^2} dx$$

The sine and cosine integrals, $si(x)$ and $Ci(x)$, are defined [24, p.878] by

(1.1)    $si(x) = -\int_x^\infty \frac{\sin t}{t} dt$

and for $x > 0$

(1.2)    $Ci(x) = -\int_x^\infty \frac{\cos t}{t} dt = \gamma + \log x + \int_0^x \frac{\cos t - 1}{t} dt$

where $\gamma$ is Euler's constant. The equivalence of these two representations is demonstrated at the end of this section.

Schlömilch [35] actually employed a slightly different sine integral $Si(x)$ which is defined in [24, p.878] and also in [1, p.231] by

$$(1.3) \qquad Si(x) = \int_0^x \frac{\sin t}{t} dt$$

We have

$$si(x) = -\int_x^\infty \frac{\sin t}{t} dt = \int_0^x \frac{\sin t}{t} dt - \int_0^\infty \frac{\sin t}{t} dt$$

and using the well-known integral from Fourier series analysis [4, p.286]

$$(1.4) \qquad \frac{\pi}{2} = \int_0^\infty \frac{\sin t}{t} dt$$

we therefore see that the two sine integrals are intimately related by

$$(1.5) \qquad si(x) = Si(x) - \frac{\pi}{2}$$

Schlömilch [35] has shown in a rather complicated manner that

$$(1.6) \qquad \int_0^\infty \frac{a \sin(tx)}{a^2 - x^2} dx = \sin(at)Ci(at) - \cos(at)Si(at)$$

$$(1.7) \qquad \int_0^\infty \frac{x \cos(tx)}{a^2 - x^2} dx = \cos(at)Ci(at) + \sin(at)Si(at)$$

Curiously, it will be noted for example that the left-hand side of (1.6) approaches zero as $a \to \infty$ whereas the right-hand side oscillates (which initially made me think that the correct integrals should involve $si(at)$ rather than $Si(at)$). We note that (1.6) appears in [7, Ch. II, (18)].

Schlömilch [35] has also shown that

$$(1.8) \qquad \int_0^\infty \frac{a \sin(tx)}{a^2 - x^2} dx = -\frac{1}{2} \pi \cos(at)$$

$$(1.9) \qquad \int_0^\infty \frac{x\cos(tx)}{a^2 - x^2}\,dx = \frac{1}{2}\pi\sin(at)$$

□

In this section we show the equivalence of the two representations for the cosine integral set out in (1.2). This analysis is based on a short paper by Omarjee [34] which in turn relies on the treatise published by Gronwall [25] in 1918.

We designate $A(x)$ as

$$A(x) = \int_0^x \frac{1-\cos t}{t}\,dt$$

and consider

$$(1.10) \qquad A(n\pi) = \int_0^{n\pi} \frac{1-\cos t}{t}\,dt$$

$$= \int_0^{1/2} \frac{1-\cos 2n\pi u}{u}\,du$$

Simple algebra gives us

$$\frac{1-\cos 2n\pi u}{u} = \pi\left[\frac{1-\cos 2n\pi u}{\sin \pi u}\right] + \left[\frac{1}{u} - \frac{\pi}{\sin \pi u}\right] - \left[\frac{1}{u} - \frac{\pi}{\sin \pi u}\right]\cos 2n\pi u$$

and we have

$$A(n\pi) = \pi\int_0^{1/2}\left[\frac{1-\cos 2n\pi u}{\sin \pi u}\right]du + \int_0^{1/2} f(u)\,du - \int_0^{1/2} f(u)\cos 2n\pi u\,du$$

where $f(u) = \frac{1}{u} - \frac{\pi}{\sin \pi u}$ and L'Hôpital's rule shows that $\lim_{u \to 0} f(u) = 0$.

We have

$$\pi\int_0^{1/2}\left[\frac{1-\cos 2n\pi u}{\sin \pi u}\right]du = 2\pi\int_0^{1/2}\sum_{k=1}^{n}\sin(2k-1)\pi u\,du$$

$$= 2\sum_{k=1}^{n} \frac{1}{2k-1}$$

and we easily determine that

$$\int_0^{1/2} f(u)\,du = \int_0^{1/2} \left[\frac{1}{u} - \frac{\pi}{\sin \pi u}\right] du$$

$$= \log \frac{u}{\tan(\pi u/2)} \bigg|_0^{1/2}$$

$$= \log \pi - 2\log 2$$

Hence we have

(1.11) $$A(n\pi) = 2\sum_{k=1}^{n} \frac{1}{2k-1} + \log \pi - 2\log 2 - \int_0^{1/2} f(u)\cos 2n\pi u\,du$$

We now consider the difference

$$c_n = A(n\pi) - \log(n\pi) = \int_0^{n\pi} \frac{1-\cos t}{t}\,dt - \log(n\pi)$$

$$= \int_0^1 \frac{1-\cos t}{t}\,dt + \int_1^{n\pi} \frac{1-\cos t}{t}\,dt - \log(n\pi)$$

$$= \int_0^1 \frac{1-\cos t}{t}\,dt - \int_1^{n\pi} \frac{\cos t}{t}\,dt$$

and we see that

$$\lim_{n\to\infty} c_n = \int_0^1 \frac{1-\cos t}{t}\,dt - \int_1^{\infty} \frac{\cos t}{t}\,dt$$

We note from (1.11) that

$$c_n = 2\sum_{k=1}^{n} \frac{1}{2k-1} - \log n - 2\log 2 - \int_0^{1/2} f(u)\cos 2n\pi u\,du$$

From the definition of Euler's constant we have

$$\lim_{n\to\infty}\left[2\sum_{k=1}^{4n}\frac{1}{k}-2\log(4n)\right]=2\gamma$$

and

$$\lim_{n\to\infty}\left[\sum_{k=1}^{2n}\frac{1}{k}-\log(2n)\right]=\gamma$$

and by subtraction we obtain

$$\lim_{n\to\infty}\left[2\sum_{k=1}^{2n}\frac{1}{2k-1}-\log n\right]=\gamma+3\log 2$$

or equivalently

(1.12) $$\lim_{n\to\infty}\left[2\sum_{k=1}^{n}\frac{1}{2k-1}-\log n\right]=\gamma+2\log 2$$

From the Riemann-Lebesgue lemma [4, p.313] we see that

(1.13) $$\lim_{n\to\infty}\int_0^{1/2} f(u)\cos 2n\pi u\, du = 0$$

which may also be obtained directly using integration by parts.

Hence we see that $\lim_{n\to\infty} c_n = \gamma$ and we obtain the known integral

(1.14) $$\gamma=\int_0^1\frac{1-\cos t}{t}dt-\int_1^\infty\frac{\cos t}{t}dt$$

which may be written as

$$=\int_0^x\frac{1-\cos t}{t}dt+\int_x^1\frac{1-\cos t}{t}dt-\int_1^x\frac{\cos t}{t}dt-\int_x^\infty\frac{\cos t}{t}dt$$

$$=\int_0^x\frac{1-\cos t}{t}dt-\log x-\int_x^\infty\frac{\cos t}{t}dt$$

and thus we have shown that for $x>0$

(1.15) $$\int_0^x \frac{1-\cos t}{t} dt = \gamma + \log x + \int_x^\infty \frac{\cos t}{t} dt$$

which proves (1.2).

□

We see from (1.10) and (1.11) that

$$\int_0^{1/2} \frac{1-\cos 2n\pi u}{u} du = 2\sum_{k=1}^n \frac{1}{2k-1} + \log \pi - 2\log 2 - \int_0^{1/2} f(u)\cos 2n\pi u\, du$$

and we have

$$\int \frac{1-\cos 2n\pi u}{u} du = \log u - Ci(2n\pi u)$$

We see from (1.2) that

$$\lim_{x \to 0}[\log x - Ci(x)] = -\gamma$$

and using

$$\log u - Ci(2n\pi u) = \log(2n\pi u) - Ci(2n\pi u) - \log(2n\pi)$$

we obtain

$$\lim_{u \to 0}[\log u - Ci(2n\pi u)] = -\gamma - \log(2n\pi)$$

Therefore we have

(1.16) $$\int_0^{1/2} \frac{1-\cos 2n\pi u}{u} du = -Ci(n\pi) + \gamma + \log(n\pi)$$

and we see that

$$\int_0^{1/2} f(u)\cos 2n\pi u\, du = 2\sum_{k=1}^n \frac{1}{2k-1} + \log \pi - 2\log 2 + Ci(n\pi) - \gamma - \log(n\pi)$$

$$= Ci(n\pi) + 2\sum_{k=1}^n \frac{1}{2k-1} - \log n - 2\log 2 - \gamma$$

Noting that $\lim_{n\to\infty}\left[2\sum_{k=1}^{n}\frac{1}{2k-1}-\log n\right]=\gamma+2\log 2$ and $\lim_{n\to\infty}Ci(n\pi)=0$, it is clear that the integral approaches zero as $n\to\infty$ (as predicated by the Riemann-Lebesgue lemma).

We note from [38, p.20] that

(1.17) $$2\sum_{k=1}^{n}\frac{1}{2k-1}=\psi\left(n+\frac{1}{2}\right)+\gamma+2\log 2$$

where $\psi(x)$ is the digamma function [sr, p.] and this gives us

(1.18) $$\int_{0}^{1/2}f(u)\cos 2n\pi u\,du=Ci(n\pi)+\psi\left(n+\frac{1}{2}\right)-\log n$$

We note from [19] that if $\phi(x)$ is of bounded variation on $[0,b]$ then we have a version of the Poisson summation formula

(1.19) $$\frac{1}{2}\phi(0+)=\int_{0}^{b}\phi(x)\,dx+2\sum_{n=1}^{\infty}\int_{0}^{b}\phi(x)\cos 2\pi nx\,dx$$

and with our definition of $f(u)$, where $f(0+)=0$, this becomes

$$\log\frac{4}{\pi}=2\sum_{n=1}^{\infty}\int_{0}^{1/2}f(u)\cos 2n\pi u\,du$$

and (1.18) gives us the series

(1.20) $$\log\frac{4}{\pi}=2\sum_{n=1}^{\infty}\left[Ci(n\pi)+\psi\left(n+\frac{1}{2}\right)-\log n\right]$$

We showed in [15] that

(1.21) $$\psi(a+1/2)=\log a+2\sum_{n=1}^{\infty}(-1)^{n}\left[\cos(2n\pi a)Ci(2n\pi a)+\sin(2n\pi a)si(2n\pi a)\right]$$

and with $a=1/2$ we get

(1.22) $$\sum_{n=1}^{\infty}Ci(n\pi)=\frac{1}{2}[\log 2-\gamma]$$

and hence we have

$$(1.23) \qquad 2\sum_{n=1}^{\infty}\left[\psi\left(n+\frac{1}{2}\right)-\log n\right]=\gamma+\log\frac{2}{\pi}$$

It may be seen from (1.12) and (1.17) that $\lim_{n\to\infty}\left[\psi\left(n+\frac{1}{2}\right)-\log n\right]$ as required for the convergence of the above series.

□

Another demonstration of (1.23) is shown below.

Merkle and Merkle [30] have shown that

$$(1.24) \qquad \sum_{n=0}^{\infty}\left[\psi(x+n)-\psi(1+n)-\frac{x-1}{1+n}\right]=(1-x)[\psi(x)+\gamma-1]$$

where a misprint in their formula has been corrected. This may be written as

$$\sum_{n=1}^{\infty}\left[\psi(x+n)-\psi(1+n)-\frac{x-1}{1+n}\right]=-x[\psi(x)+\gamma]+2(x-1)$$

With $x=1/2$ we get

$$\sum_{n=1}^{\infty}\left[\psi\left(n+\frac{1}{2}\right)-\psi(1+n)+\frac{1}{2(1+n)}\right]=-\frac{1}{2}\left[\psi\left(\frac{1}{2}\right)+\gamma+2\right]$$

which may be written as

$$\sum_{n=1}^{\infty}\left[\psi\left(n+\frac{1}{2}\right)-\log n+\log n-\psi(1+n)+\frac{1}{2(1+n)}\right]=\log 2-1$$

Equivalently we have

$$\sum_{n=1}^{\infty}\left[\psi\left(n+\frac{1}{2}\right)-\log n\right]-\sum_{n=1}^{\infty}\left[\psi(1+n)-\log n-\frac{1}{2(1+n)}\right]=\log 2-1$$

It is shown in [19] that

$$(1.25) \qquad \sum_{n=1}^{\infty}\left[\psi(1+n)-\log n-\frac{1}{2(1+n)}\right]=1+\frac{1}{2}\gamma-\frac{1}{2}\log(2\pi)$$

and we again deduce (1.23)

$$2\sum_{n=1}^{\infty}\left[\psi\left(n+\frac{1}{2}\right)-\log n\right]=\gamma+\log\frac{2}{\pi}$$

□

Referring to (1.18) we make the summation

$$S=\sum_{n=1}^{\infty}\frac{1}{n^2}\int_0^{1/2}f(u)\cos 2n\pi u\,du=\sum_{n=1}^{\infty}\frac{Ci(n\pi)}{n^2}+\sum_{n=1}^{\infty}\frac{1}{n^2}\psi\left(n+\frac{1}{2}\right)-\sum_{n=1}^{\infty}\frac{\log n}{n^2}$$

We showed in equation (6.130v) in [15] that

(1.26) $$\sum_{n=1}^{\infty}\frac{Ci(n\pi)}{n^2}=-\frac{1}{6}\pi^2\log 2-\frac{5}{24}\pi^2-2\pi^2\varsigma'(-1)$$

and de Doelder [20] showed that

(1.27) $$\sum_{n=1}^{\infty}\frac{1}{n^2}\psi\left(n+\frac{1}{2}\right)=\frac{7}{2}\varsigma(3)-(\gamma+2\log 2)\varsigma(2)$$

and this identity was subsequently generalised by Coffey [13].

This gives us

$$S=-\frac{1}{6}\pi^2\log 2-\frac{5}{24}\pi^2-2\pi^2\varsigma'(-1)+\frac{7}{2}\varsigma(3)-(\gamma+2\log 2)\varsigma(2)+\varsigma'(2)$$

and using

$$\varsigma'(-1)=\frac{1}{12}(1-\gamma-\log 2\pi)+\frac{1}{2\pi^2}\varsigma'(2)$$

this becomes

$$S=-\frac{3}{8}\pi^2+\frac{7}{2}\varsigma(3)+\varsigma(2)(\log\pi-2\log 2)$$

Assuming that interchanging the order of summation and integration is valid, we have

$$\sum_{n=1}^{\infty} \frac{1}{n^2} \int_0^{1/2} f(u) \cos 2n\pi u \, du = \int_0^{1/2} f(u) \sum_{n=1}^{\infty} \frac{\cos 2n\pi u}{n^2} \, du$$

and using the well known Fourier series [4, p.338]

$$\sum_{n=1}^{\infty} \frac{\cos 2n\pi u}{n^2} = \pi^2 B_2(u) = \pi^2 \left( u^2 - u + \frac{1}{6} \right)$$

we obtain

$$\sum_{n=1}^{\infty} \frac{1}{n^2} \int_0^{1/2} f(u) \cos 2n\pi u \, du = \pi^2 \int_0^{1/2} f(u) \left( u^2 - u + \frac{1}{6} \right) du$$

We see that

$$\pi^2 \int_0^{1/2} f(u) \left( u^2 - u + \frac{1}{6} \right) du = \pi^2 \int_0^{1/2} \left( \frac{1}{u} - \frac{\pi}{\sin \pi u} \right) \left( u^2 - u + \frac{1}{6} \right) du$$

$$= \frac{1}{6} \pi^2 \int_0^{1/2} \left( \frac{1}{u} - \frac{\pi}{\sin \pi u} \right) du + \pi^2 \int_0^{1/2} \frac{\pi u^2}{\sin \pi u} du - \pi^2 \int_0^{1/2} \frac{\pi u}{\sin \pi u} du$$

$$+ \pi^2 \int_0^{1/2} (u - 1) \, du$$

We showed in [15] that

(1.28) $$\int_0^{\pi/2} \frac{x}{\sin x} dx = 2 \sum_{n=0}^{\infty} \frac{(-1)^n}{(2n+1)^2} = 2G$$

where $G$ is Catalan's constant. This result is well-known and an alternative proof is contained, for example, in Bradley's website [10].

We also showed that [15]

(1.29) $$\int_0^{\pi/2} \frac{x^2}{\sin x} dx = 2\pi G - \frac{7}{2} \varsigma(3)$$

The evaluation of this integral is also contained in Bradley's website [10] and an alternative proof was provided by De Doelder [20].

The desired equivalence follows by simple algebra.

In the following section we show how the sine and cosine integrals may be employed to derive Kummer's Fourier series for $\log \Gamma(x)$.

## 2. A new derivation of Kummer's Fourier series for $\log \Gamma(x)$

We showed in [18] that for $a > 0$ and $n \geq 1$

$$(2.1) \quad \int_0^1 \log \Gamma(x+a) \sin 2n\pi x \, dx = -\frac{1}{2n\pi}\left[\log a - \cos(2n\pi a) Ci(2n\pi a) - \sin(2n\pi a) si(2n\pi a)\right]$$

$$(2.2) \quad \int_0^1 \log \Gamma(x+a) \cos 2n\pi x \, dx = -\frac{1}{2n\pi}\left[-\sin(2n\pi a) Ci(2n\pi a) + \cos(2n\pi a) si(2n\pi a)\right]$$

which correct the entries reported in [24, p.650, 6.443.5 & 6.443.6]. As noted in [18], both of these integrals also hold in the limit as $a \to 0$.

We designate $F(a)$ by

$$F(a) = \int_0^1 \log \Gamma(a+x) \, dx$$

and differentiation gives us

$$F'(a) = \int_0^1 \psi(a+x) \, dx = \log \Gamma(a+1) - \log \Gamma(a) = \log a$$

Integration then results in

$$F(u) - F(0) = u \log u - u$$

and using Raabe's integral [38, p.18]

$$F(0) = \int_0^1 \log \Gamma(x) \, dx = \frac{1}{2}\log(2\pi)$$

we see that

$$\int_0^1 \log \Gamma(a+x) \, dx = \frac{1}{2}\log(2\pi) + a \log a - a$$

as reported in [38, p.207].

We then deduce the Fourier series for $\log \Gamma(a+x)$ which is valid for $0 < x < 1$

(2.3) $\quad \log \Gamma(a+x) = \frac{1}{2}\log(2\pi) + a\log a - a$

$$-\sum_{n=1}^{\infty} \frac{1}{n\pi}\left[-\sin(2n\pi a)Ci(2n\pi a) + \cos(2n\pi a)si(2n\pi a)\right]\cos(2n\pi x)$$

$$-\sum_{n=1}^{\infty} \frac{1}{n\pi}\left[\log a - \cos(2n\pi a)Ci(2n\pi a) - \sin(2n\pi a)si(2n\pi a)\right]\sin(2n\pi x)$$

Letting $a \to a/2$ and $x \to a/2$ in (2.3) gives us

$$\log \Gamma(a) = \frac{1}{2}\log(2\pi) + \frac{1}{2}a\log(a/2) - \frac{1}{2}a$$

$$-\sum_{n=1}^{\infty} \frac{1}{n\pi}\left[\log(a/2) - 2\sin(n\pi a)\cos(n\pi a)Ci(n\pi a) + \left\{\cos^2(n\pi a) - \sin^2(n\pi a)\right\}si(n\pi a)\right]$$

Equation (2.3) also applies in the limit as $a \to 0$ because referring to (1.2) we have

$$\lim_{y \to 0}[\cos y\, Ci(y) - \log y] = \lim_{y \to 0}\left[\gamma \cos y + \log y[\cos y - 1] + \cos y \int_0^y \frac{\cos t - 1}{t}dt\right]$$

and, applying L'Hôpital's rule, we see that

$$\lim_{y \to 0}[\log y(\cos y - 1)] = \lim_{y \to 0}\left[y \log y \frac{\cos y - 1}{y}\right]$$

$$= \lim_{y \to 0}[y \log y]\lim_{y \to 0}\left[\frac{\cos y - 1}{y}\right]$$

$$= -\lim_{y \to 0}[y \log y]\lim_{y \to 0}[\sin y] = 0$$

We therefore obtain

(2.4) $\quad \lim_{y \to 0}[\cos y\, Ci(y) - \log y] = \gamma$

or equivalently

$$\lim_{a \to 0}[\cos(2n\pi a) Ci(2n\pi a) - \log(2n\pi a)] = \gamma$$

which results in

(2.5) $$\lim_{a \to 0}[\cos(2n\pi a) Ci(2n\pi a) - \log a] = \gamma + \log(2n\pi)$$

Similarly we have

$$\lim_{y \to 0}[\sin y\, Ci(y)] = \lim_{y \to 0}\left[\gamma \sin y + \sin y \log y + \sin y \int_0^y \frac{\cos t - 1}{t} dt\right]$$

and, since $\lim_{y \to 0} \sin y \log y = \lim_{y \to 0} \frac{\sin y}{y} \lim_{y \to 0} y \log y = 0$, this results in

(2.6) $$\lim_{a \to 0} \sin(2n\pi a) Ci(2n\pi a) = 0$$

Therefore we have the limit as $a \to 0$

$$\log \Gamma(x) = \frac{1}{2} \log(2\pi) - \sum_{n=1}^{\infty} \frac{si(0) \cos(2n\pi x)}{n\pi} + \sum_{n=1}^{\infty} \frac{\gamma + \log(2n\pi)}{n\pi} \sin(2n\pi x)$$

and using (1.1) $si(0) = -\pi/2$ we obtain Kummer's Fourier series for $\log \Gamma(x)$ ([27] and [38, p.17]) which is valid for $0 < x < 1$

(2.7) $$\log \Gamma(x) = \frac{1}{2} \log(2\pi) + \sum_{n=1}^{\infty}\left[\frac{1}{2n} \cos(2n\pi x) + \frac{\gamma + \log(2n\pi)}{n\pi} \sin(2n\pi x)\right]$$

With $a = 1$ in (2.3) we obtain Nielsen's representation [31, p.79] which is valid for $0 < x < 1$

(2.8) $$\log \Gamma(1+x) = \frac{1}{2} \log(2\pi) - 1 + \sum_{n=1}^{\infty} \frac{Ci(2n\pi) \sin(2n\pi x) - si(2n\pi) \cos(2n\pi x)}{n\pi}$$

Note that in this case, the sine and cosine integral functions do not contain the variable $x$ in their arguments.

Equation (2.3) may be written in a more compact form as

$$\log \Gamma(a+x) = \frac{1}{2} \log(2\pi) + a \log a - a - \log a \sum_{n=1}^{\infty} \frac{\sin(2n\pi x)}{n\pi}$$

$$+\sum_{n=1}^{\infty}\frac{Ci(2n\pi a)\sin 2n\pi(a+x)-si(2n\pi a)\cos 2n\pi(a+x)}{n\pi}$$

Using the familiar Fourier series [12, p.241] which is valid for $0 < x < 2\pi$

$$\frac{1}{2}(\pi-x)=\sum_{n=1}^{\infty}\frac{\sin nx}{n}$$

we obtain for $0 < x < 1$

(2.9) $\quad\log\Gamma(a+x)=\frac{1}{2}\log(2\pi)+\left(a+x-\frac{1}{2}\right)\log a-a$

$$+\sum_{n=1}^{\infty}\frac{\sin[2n\pi(a+x)]Ci(2n\pi a)-\cos[2n\pi(a+x)]si(2n\pi a)}{n\pi}$$

and with $a=1$ we obtain (2.8) again.

By symmetry we have for $0 < a < 1$

(2.10) $\quad\log\Gamma(a+x)=\frac{1}{2}\log(2\pi)+\left(a+x-\frac{1}{2}\right)\log x-x$

$$+\sum_{n=1}^{\infty}\frac{\sin[2n\pi(a+x)]Ci(2n\pi x)-\cos[2n\pi(a+x)]si(2n\pi x)}{n\pi}$$

With $x=0$ in (2.9) we obtain for $0 < a < 1$

(2.11) $\quad\log\Gamma(a)=\frac{1}{2}\log(2\pi)+\left(a-\frac{1}{2}\right)\log a-a$

$$+\sum_{n=1}^{\infty}\frac{\sin(2n\pi a)Ci(2n\pi a)-\cos(2n\pi a)si(2n\pi a)}{n\pi}$$

A different derivation of (2.11) is given in [18] using Bourguet's formula [40, p.261]

$$\log\Gamma(a)=\frac{1}{2}\log(2\pi)+\left(a-\frac{1}{2}\right)\log a-a+\frac{1}{\pi}\sum_{n=1}^{\infty}\int_{0}^{\infty}\frac{\sin(2n\pi x)}{n(x+a)}dx$$

which may be derived using the Euler-Maclaurin summation formula (see for example Knopp's book [26, p.530]). The formula (2.11) was also given by Nörlund [33, p.114].

Letting $a = 1/2$ in (2.10) results in

(2.12)

$$\log \Gamma\left(x + \frac{1}{2}\right) = \frac{1}{2}\log(2\pi) + x \log x - x + \sum_{n=1}^{\infty} \frac{(-1)^n}{n\pi}[\sin 2n\pi x \, Ci(2n\pi x) - \cos 2n\pi x \, si(2n\pi x)]$$

Letting $x \to 1 - x$ in (2.7) results in

$$\log \Gamma(1-x) = \frac{1}{2}\log(2\pi) + \sum_{n=1}^{\infty} \left[\frac{1}{2n}\cos(2n\pi x) - \frac{\gamma + \log(2n\pi)}{n\pi}\sin(2n\pi x)\right]$$

so that

$$\log \Gamma(x) + \log \Gamma(1-x) = \log(2\pi) + \sum_{n=1}^{\infty} \frac{\cos(2n\pi x)}{n}$$

and, as pointed out by Berndt [8], using (3.4.2) we easily obtain Euler's reflection formula for the gamma function

$$\Gamma(x)\Gamma(1-x) = \frac{\pi}{\sin \pi x}$$

We also have the difference

$$\log \Gamma(x) - \log \Gamma(1-x) = 2\sum_{n=1}^{\infty} \frac{\gamma + \log(2n\pi)}{n\pi}\sin(2n\pi x)$$

□

Largely for my own benefit, we set out a heuristic exposition of Parseval's Theorem in the case where $[0,1]$ is the interval of integration:

Letting

$$f(x) = a_0 + \sum_{n=1}^{\infty}(a_n \cos 2\pi nx + b_n \sin 2\pi nx)$$

we see that

$$\int_0^1 f(x)\,dx = a_0$$

$$\int_0^1 f(x)\cos 2\pi nx\,dx = a_n \int_0^1 \cos^2 2\pi nx\,dx = \frac{1}{2}a_n$$

$$\int_0^1 f(x)\sin 2\pi nx\,dx = b_n \int_0^1 \sin^2 2\pi nx\,dx = \frac{1}{2}b_n$$

and thus we have

$$\int_0^1 f^2(x)\,dx = a_0^2 + \frac{1}{2}\sum_{n=1}^{\infty}(a_n^2 + b_n^2)$$

or

$$\int_0^1 f^2(x)\,dx = \left[\int_0^1 f(x)\,dx\right]^2 + 2\sum_{n=1}^{\infty}\left(\left[\int_0^1 f(x)\cos 2\pi nx\,dx\right]^2 + \left[\int_0^1 f(x)\sin 2\pi nx\,dx\right]^2\right)$$

$\square$

Applying Parseval's theorem to (2.3) gives us

$$\int_0^1 \log^2 \Gamma(x+a)\,dx = \left[\frac{1}{2}\log(2\pi) + a\log a - a\right]^2 + \frac{1}{2\pi^2}\sum_{n=1}^{\infty}\frac{f_n(a)}{n^2}$$

where

$$f_n(a) = Ci^2(2n\pi a) + si^2(2n\pi a) + \log^2 a - 2[\cos(2n\pi a)Ci(2n\pi a) + \sin(2n\pi a)si(2n\pi a)]\log a$$

We see that $f_n(1) = Ci^2(2n\pi) + si^2(2n\pi)$ and thus we have

$$\int_0^1 \log^2 \Gamma(x+1)\,dx = \left[\frac{1}{2}\log(2\pi) - 1\right]^2 + \frac{1}{2\pi^2}\sum_{n=1}^{\infty}\frac{Ci^2(2n\pi) + si^2(2n\pi)}{n^2}$$

With $a = 1/2$ we get

$$\int_0^1 \log^2 \Gamma\left(x+\frac{1}{2}\right)dx = \frac{1}{4}[\log \pi - 1]^2 + \frac{1}{2\pi^2}\sum_{n=1}^{\infty}\frac{Ci^2(n\pi) + si^2(n\pi)}{n^2}$$

$$+ \frac{1}{12}\log^2 2 + \frac{\log 2}{\pi^2}\sum_{n=1}^{\infty}\frac{(-1)^n Ci(n\pi)}{n^2}$$

With $a = 0$ we have

$$f_n(0) = \frac{\pi^2}{4} + \lim_{a \to 0}[Ci^2(2n\pi a) + \log^2 a - 2\cos(2n\pi a)Ci(2n\pi a)\log a]$$

It is shown in [18] that

$$\lim_{a \to 0}[Ci(\mu a) - \cos(\mu a)\log a] = \gamma + \log \mu$$

Therefore we see that

$$\lim_{a \to 0}[Ci(\mu a) - \cos(\mu a)\log a]^2 = [\gamma + \log \mu]^2$$

and we have

$$[Ci(\mu a) - \cos(\mu a)\log a]^2 = Ci^2(\mu a) + \cos^2(\mu a)\log^2 a - 2\cos(\mu a)Ci(\mu a)\log a$$

$$= Ci^2(\mu a) + \log^2 a - 2\cos(\mu a)Ci(\mu a)\log a$$

$$+ [\cos^2(\mu a) - 1]\log^2 a$$

We see that

$$[\cos^2(\mu a) - 1]\log^2 a = [\cos(\mu a) + 1][\cos(\mu a) - 1]\log^2 a$$

$$= [\cos(\mu a) + 1]\frac{\cos(\mu a) - 1}{a^2}(a\log a)^2$$

and using L'Hôpital's rule we see that

$$\lim_{a \to 0}[\cos^2(\mu a) - 1]\log^2 a = 0$$

Therefore we obtain

$$f_n(0) = \frac{\pi^2}{4} + [\gamma + \log(2n\pi)]^2$$

which results in

$$\int_0^1 \log^2 \Gamma(x)\,dx = \frac{1}{4}\log^2(2\pi) + \frac{\pi^2}{48} + \frac{1}{2\pi^2}\sum_{n=1}^{\infty}\frac{1}{n^2}\left(\gamma^2 + 2\gamma\log(2n\pi) + \log^2(2n\pi)\right)$$

Noting that

$$\gamma^2 + 2\gamma \log(2n\pi) + \log^2(2n\pi) = \gamma^2 + 2\gamma \log(2\pi) + 2[\gamma + \log(2\pi)]\log n + \log^2 n + \log^2(2\pi)$$

we immediately deduce that

$$\int_0^1 \log^2 \Gamma(x)\, dx = \frac{\gamma^2}{12} + \frac{\pi^2}{48} + \frac{1}{6}\gamma \log(2\pi) + \frac{1}{3}\log^2(2\pi) - [\gamma + \log(2\pi)]\frac{\varsigma'(2)}{\pi^2} + \frac{\varsigma''(2)}{2\pi^2}$$

which was previously obtained by Espinosa and Moll [22].

We shall see in (7.5) that

$$\frac{1}{\pi^2} \sum_{n=1}^{\infty} \frac{\cos(2n\pi a)\, Ci(2n\pi a) + \sin(2n\pi a)\, si(2n\pi a)}{n^2}$$

$$= 2[\log G(1+a) - a \log \Gamma(a)] + \left[a(a-1) + \frac{1}{6}\right]\log a - \frac{1}{2}a^2 + 2\left[\frac{1}{12} - \varsigma'(-1)\right]$$

and hence we obtain

$$\int_0^1 \log^2 \Gamma(x+a)\, dx = \left[\frac{1}{2}\log(2\pi) + a \log a - a\right]^2 + \frac{1}{2\pi^2}\sum_{n=1}^{\infty} \frac{Ci^2(2n\pi a) + si^2(2n\pi a)}{n^2} + \frac{1}{12}\log a$$

$$- \left\{ 2[\log G(1+a) - a \log \Gamma(a)] + \left[a(a-1) + \frac{1}{6}\right]\log a - \frac{1}{2}a^2 + 2\left[\frac{1}{12} - \varsigma'(-1)\right]\right\}\log a$$

A different method of deriving Kummer's Fourier series for $\log \Gamma(x)$ is shown below.

**3. Fourier coefficients for $\log x$**

The following exposition is adopted from Nielsen's book [31] with a little added embellishment.

We have the Fourier series for $\log x$

$$\log x = \frac{1}{2}a_0 + \sum_{n=1}^{\infty}(a_n \cos nx + b_n \sin nx)$$

where, as proved below, we have

(3.1) $$a_n = \frac{1}{\pi}\int_0^{2\pi} \log x \cdot \cos nx \, dx = -\frac{1}{n\pi}\left[si(2n\pi) + \frac{\pi}{2}\right] = -\frac{Si(2n\pi)}{n\pi}$$

(3.2) $$b_n = \frac{1}{\pi}\int_0^{2\pi} \log x \cdot \sin nx \, dx = \frac{1}{n\pi}\left[Ci(2n\pi) - \gamma - \log(2n\pi)\right]$$

First of all we consider the integral

$$\int_0^x \log t \cdot \cos nt \, dt = \log t \frac{\sin nt}{n}\bigg|_0^x - \int_0^x \frac{\sin nt}{nt} dt$$

We have

$$\lim_{t \to 0} \log t \sin nt = \lim_{t \to 0} t \log t \frac{\sin nt}{t} = \lim_{t \to 0} t \log t \lim_{t \to 0} \frac{\sin nt}{t} = 0$$

Therefore we get

$$\int_0^x \log t \cdot \cos nt \, dt = \frac{\sin nx \log x}{n} - \int_0^x \frac{\sin nt}{nt} dt$$

$$= \frac{\sin nx \log x}{n} - \frac{1}{n}\int_0^{nx} \frac{\sin u}{u} du$$

and hence we obtain

$$\int_0^x \log t \cdot \cos nt \, dt = \frac{\sin nx \log x}{n} - \frac{Si(nx)}{n}$$

With integration by parts we get

$$\int_\alpha^x \log t \cdot \sin at \, dt = -\log t \frac{\cos at}{a}\bigg|_\alpha^x + \int_\alpha^x \frac{\cos at}{at} dt$$

$$= -\log t \frac{\cos at}{a}\bigg|_\alpha^x + \int_\alpha^x \frac{\cos at - 1}{at} dt + \int_\alpha^x \frac{1}{at} dt$$

$$= -\log x \frac{\cos ax}{a} + \int_\alpha^{ax} \frac{\cos u - 1}{u} du + \frac{\log \alpha}{a}(\cos a\alpha - 1) + \frac{\log x}{a}$$

Therefore, in the limit as $\alpha \to 0$, we get

$$\int_0^x \log t \cdot \sin at \, dt = -\log x \frac{\cos ax}{a} + \int_0^{ax} \frac{\cos u - 1}{u} du + \frac{\log x}{a}$$

and using (1.15) we have

(3.3) $$\int_0^x \log t \cdot \sin at \, dt = \frac{1}{a}\left[Ci(ax) - \gamma - \cos(ax)\log x - \log a\right]$$

as reported in Nielsen's book [31, p.12],

We therefore obtain

$$b_n = \frac{1}{\pi} \int_0^{2\pi} \log x \cdot \sin nx \, dx = \frac{1}{n\pi}\left[Ci(2n\pi) - \gamma - \log(2n\pi)\right]$$

We also have

$$\frac{1}{2}a_0 = \frac{1}{2\pi}\int_0^{2\pi} \log x \, dx = \log(2\pi) - 1$$

and we obtain the Fourier series (where we have made the substitution $x \to 2\pi x$)

$$\log(2\pi x) = \log(2\pi) - 1 + \sum_{n=1}^{\infty} \frac{[Ci(2n\pi) - \gamma - \log(2n\pi)]\sin(2n\pi x) - [si(2n\pi) + \pi/2]\cos(2n\pi x)}{n\pi}$$

so that for $0 < x < 1$

(3.4) $$\log x = -1 + \sum_{n=1}^{\infty} \frac{[Ci(2n\pi) - \gamma - \log(2n\pi)]\sin(2n\pi x) - [si(2n\pi) + \pi/2]\cos(2n\pi x)}{n\pi}$$

Making use of the familiar Fourier series [12, p.241]

(3.4.1) $$\frac{1}{2}(\pi - x) = \sum_{n=1}^{\infty} \frac{\sin nx}{n} \qquad (0 < x < 2\pi)$$

(3.4.2) $$\log\left[2\sin(x/2)\right] = -\sum_{n=1}^{\infty} \frac{\cos nx}{n} \qquad (0 < x < 2\pi)$$

we therefore obtain the Fourier series

$$(3.5) \quad 1+\log x - \frac{1}{2}\log(2\sin \pi x)+(\gamma+\log 2\pi)\left(\frac{1}{2}-x\right)$$

$$= \sum_{n=1}^{\infty} \frac{[Ci(2n\pi)-\log n]\sin 2n\pi x - si(2n\pi)\cos 2n\pi x}{n\pi}$$

In (3.5) we let $x \to 1-x$ and combine the resulting identity with (3.5) to obtain

$$(3.6) \quad \log x + \log(1-x) - \log(2\sin \pi x) + 2 = -\frac{2}{\pi}\sum_{n=1}^{\infty} \frac{si(2n\pi)\cos 2n\pi x}{n}$$

and subtraction results in

$$(3.7) \quad \log x - \log(1-x) + (\gamma + \log 2\pi)(1-2x) = \frac{2}{\pi}\sum_{n=1}^{\infty} \frac{[Ci(2n\pi)-\log n]\sin 2n\pi x}{n}$$

Combining (3.4) with (2.8)

$$\log \Gamma(x) + \log x = \frac{1}{2}\log(2\pi) - 1 + \sum_{n=1}^{\infty} \frac{\sin(2n\pi x)Ci(2n\pi) - \cos(2n\pi x)si(2n\pi)}{n\pi}$$

we obtain another derivation of Kummer's Fourier series for $\log \Gamma(x)$

$$\log \Gamma(x) = \frac{1}{2}\log(2\pi) + \sum_{n=1}^{\infty}\left[\frac{1}{2n}\cos(2n\pi x) + \frac{\gamma + \log(2n\pi)}{n\pi}\sin(2n\pi x)\right]$$

Integrating (3.4) gives us

$$x\log x = 2\sum_{n=1}^{\infty} \frac{[Ci(2n\pi)-\gamma-\log(2n\pi)][1-\cos(2n\pi x)] - [si(2n\pi)+\pi/2]\sin(2n\pi x)}{(2n\pi)^2}$$

and with $x = 1/2$ we have

$$\log 2 = 4\sum_{n=1}^{\infty} \frac{[Ci(2n\pi)-\gamma-\log(2n\pi)][(-1)^n - 1]}{(2n\pi)^2}$$

**4. The digamma function**

We showed in [18] that for $a > 0$

$$(4.1) \quad \int_0^1 \psi(x+a)\cos 2n\pi x\, dx = \sin(2n\pi a)si(2n\pi a) + \cos(2n\pi a)Ci(2n\pi a)$$

and for $a \geq 0$

(4.2) $\quad \int_0^1 \psi(x+a)\sin 2n\pi x\, dx = -\sin(2n\pi a)Ci(2n\pi a) + \cos(2n\pi a)si(2n\pi a)$

which correct the entries in [24, p.652, 6.467.1 & 2]. These integrals may also be obtained directly by differentiating (2.1) and (2.2).

Prima facie, we would therefore expect to have the Fourier series for $0 < x < 1$

(4.2.1) $\quad \psi(x+a) = \log a + 2\sum_{n=1}^{\infty}[\sin(2n\pi a)si(2n\pi a) + \cos(2n\pi a)Ci(2n\pi a)]\cos 2n\pi x$

$$-2\sum_{n=1}^{\infty}[\sin(2n\pi a)Ci(2n\pi a) - \cos(2n\pi a)si(2n\pi a)]\sin 2n\pi x$$

As noted in [19], using the Poisson summation formula, we obtain

(4.3) $\quad \psi(a) = \log a - \dfrac{1}{2a} + 2\sum_{n=1}^{\infty}[\cos(2n\pi a)Ci(2n\pi a) + \sin(2n\pi a)si(2n\pi a)]$

and this may be derived by differentiating (2.11). Equation (4.3) concurs with Nörlund's analysis [33, p.108]. It should be noted that this is not a Fourier series.

Applying Parseval's theorem to (4.2.1) gives us

(4.3.1) $\quad \int_0^1 \psi^2(x+a)\, dx = \log^2 a + 2\sum_{n=1}^{\infty}[Ci^2(2n\pi a) + si^2(2n\pi a)]$

so that

(4.3.2) $\quad \int_0^1 \psi^2(x+1)\, dx = 2\sum_{n=1}^{\infty}[Ci^2(2n\pi) + si^2(2n\pi)]$

$\square$

Lerch [28] proved the following theorem in 1895:

Let

$$f(x) = \sum_{n=1}^{\infty} \frac{c_n}{n}\sin 2\pi nx$$

be convergent for $0 < x < 1$; then the derivative of $f(x)$ is given in this interval by

(4.4) $$f'(x) \cdot \frac{\sin \pi x}{\pi} = \sum_{n=0}^{\infty} (c_n - c_{n+1}) \sin(2n+1)\pi x$$

where $c_0 = 0$ provided that the last series converges uniformly for $\varepsilon \leq x \leq 1-\varepsilon$ for all $\varepsilon > 0$. Equation (4.4) may be written as

$$f'(x) \cdot \frac{\sin \pi x}{\pi} = -c_1 \sin \pi x + \sum_{n=1}^{\infty} (c_n - c_{n+1}) \sin(2n+1)\pi x$$

Subject to the same conditions, we also have

$$g(x) = \sum_{n=1}^{\infty} \frac{c_n}{n} \cos 2\pi n x$$

then

$$g'(x) \cdot \frac{\sin \pi x}{\pi} = \sum_{n=0}^{\infty} (c_n - c_{n+1}) \cos(2n+1)\pi x$$

or

(4.5) $$g'(x) \cdot \frac{\sin \pi x}{\pi} = -c_1 \cos \pi x + \sum_{n=1}^{\infty} (c_n - c_{n+1}) \cos(2n+1)\pi x$$

Applying (4.4) to (3.6) results in

(4.6) $$\left[ \frac{1}{x} + \frac{1}{1-x} - 2(\gamma + \log 2\pi) \right] \sin \pi x$$

$$= -2Ci(2\pi) + 2\sum_{n=1}^{\infty} \left[ Ci(2n\pi) - Ci(2[n+1]\pi) + \log \frac{n+1}{n} \right] \sin(2n+1)\pi x$$

provided that the last series converge uniformly for $\varepsilon \leq x \leq 1-\varepsilon$ for all $\varepsilon > 0$.

With $x = 1/2$ we obtain

$$2 - (\gamma + \log 2\pi) = -Ci(2\pi) + \sum_{n=1}^{\infty} (-1)^n \left[ Ci(2n\pi) - Ci(2[n+1]\pi) + \log \frac{n+1}{n} \right]$$

$$= -Ci(2\pi) + \sum_{n=1}^{\infty} (-1)^n \left[ Ci(2n\pi) - Ci(2[n+1]\pi) \right] - \sum_{n=1}^{\infty} (-1)^{n+1} \log \frac{n+1}{n}$$

We know from [37] that

$$\log\frac{\pi}{2} = \sum_{n=1}^{\infty}(-1)^{n+1}\log\frac{n+1}{n}$$

and we have

$$\sum_{n=1}^{\infty}(-1)^{n}\left[Ci(2n\pi) - Ci(2[n+1]\pi)\right] = Ci(2\pi) + 2\sum_{n=1}^{\infty}(-1)^{n}Ci(2n\pi)$$

This results in

(4.7) $$\sum_{n=1}^{\infty}(-1)^{n}Ci(2n\pi) = 1 - \frac{1}{2}\gamma - \log 2$$

which concurs with Nielsen's book [31, p.80].

The following related result was obtained in [trig]

(4.8) $$\sum_{n=1}^{\infty}(-1)^{n}Ci(n\pi) = \frac{1}{2}\left(1 - \gamma - \log 2\right)$$

Applying (4.5) to (3.5) results in

$$\frac{1}{2}\left[\frac{1}{x} - \frac{1}{1-x} - \pi\cot\pi x\right]\sin\pi x = si(2\pi)\cos\pi x - \sum_{n=1}^{\infty}[si(2n\pi) - si(2[n+1]\pi)]\cos(2n+1)\pi x$$

□

Kummer's formula (2.7) may be written as

(4.9) $$\log\Gamma(x) = \frac{1}{2}\log\pi - \frac{1}{2}\log\sin\pi x + \frac{1}{2}(\gamma + \log 2\pi)(1-2x) + \sum_{n=1}^{\infty}\frac{\log n \sin 2\pi nx}{\pi n}$$

and applying Lerch's differentiation formula (4.4) we obtain Lerch's trigonometric series expansion for the digamma function for $0 < x < 1$ (see for example Gronwall's paper [25, p.105] and Nielsen's book, Die Gammafunktion, [32, p.204])

(4.10) $$\psi(x)\sin\pi x + \frac{\pi}{2}\cos\pi x + (\gamma + \log 2\pi)\sin\pi x = -\sum_{n=1}^{\infty}\sin(2n+1)\pi x \cdot \log\frac{n+1}{n}$$

**5. Miscellaneous integrals**

Abramowitz and Stegun [1, p.232] define the auxiliary functions $f(x)$ and $g(x)$ as

(5.1) $$f(x) = -\cos x \, si(x) + \sin x \, Ci(x) = \int_0^\infty \frac{\sin y}{y+x} dy$$

(5.2) $$g(x) = -\cos x \, Ci(x) - \sin x \, si(x) = \int_0^\infty \frac{\cos y}{y+x} dy$$

and report that for $\text{Re}(x) > 0$

(5.3) $$f(x) = \int_0^\infty \frac{e^{-xu}}{1+u^2} du$$

(5.4) $$g(x) = \int_0^\infty \frac{u e^{-xu}}{1+u^2} du$$

It is easily seen that $f'(x) = -g(x)$.

The above results may be easily derived by considering the double integral

$$I = \int_0^\infty \int_0^\infty e^{-(a+y)x} \sin y \, dx \, dy$$

where integrating with respect to $x$ gives us

$$\int_0^\infty e^{-(a+y)x} dx = \frac{1}{a+y}$$

and thus we have

$$I = \int_0^\infty \frac{\sin y}{a+y} dy$$

Similarly, integrating with respect to $y$ gives us

$$\int_0^\infty e^{-(a+y)x} \sin y \, dy = \frac{e^{-ax}}{2i} \int_0^\infty e^{-yx}(e^{iy} - e^{-iy}) dy = \frac{e^{-ax}}{1+x^2}$$

Therefore we see that

(5.5) $$\int_0^\infty \frac{e^{-ax}}{1+x^2} dx = \int_0^\infty \frac{\sin y}{a+y} dy$$

and the validity of the operation

$$\int_0^\infty dx \int_0^\infty e^{-(a+y)x} \sin y \, dy = \int_0^\infty dy \int_0^\infty e^{-(a+y)x} \sin y \, dx$$

is confirmed by [6, p.282].

The formula

(5.6) $$\int_0^\infty \frac{xe^{-ax}}{1+x^2} dx = \int_0^\infty \frac{\cos y}{a+y} dy$$

may be derived in a similar fashion.

It may be seen that

$$f(ax) = -\cos(ax)\,si(ax) + \sin(ax)\,Ci(ax) = \int_0^\infty \frac{\sin y}{y+ax} dy$$

and we have

$$\int_0^\infty \frac{\sin y}{y+ax} dy = \int_0^\infty \frac{\sin(av)}{v+x} dv$$

so that

$$f(ax) = -\cos(ax)\,si(ax) + \sin(ax)\,Ci(ax) = \int_0^\infty \frac{\sin(av)}{v+x} dv$$

Similarly we have

$$g(ax) = -[\cos(ax)\,Ci(ax) + \sin(ax)\,si(ax)] = \int_0^\infty \frac{\cos(av)}{v+x} dv$$

We see from (5.5) that

$$f(ax) = \int_0^\infty \frac{e^{-axu}}{1+u^2} du$$

$$= \int_0^\infty \frac{e^{-xv}}{1+v^2/a^2} dv$$

$$= a \int_0^\infty \frac{e^{-xv}}{a^2 + v^2} dv$$

We accordingly obtain the known integral [24, p.338, 3.354.1]

(5.7) $$\int_0^\infty \frac{e^{-xv}}{a^2 + v^2} dv = \frac{1}{a}[\sin(ax)Ci(ax) - \cos(ax)si(ax)]$$

which was originally derived in another paper by Schlömilch [36] in 1846 (his two papers [35] and [36] are actually printed in same the journal alongside each other).

Similarly we find that

(5.8) $$g(ax) = \int_0^\infty \frac{ve^{-xv}}{a^2 + v^2} dv$$

so that for $a, \mu > 0$ we have [24, p.338, 3.354.2]

(5.9) $$\int_0^\infty \frac{ve^{-\mu v}}{a^2 + v^2} dv = -[\cos(a\mu)Ci(a\mu) + \sin(a\mu)si(a\mu)]$$

We see that (5.7) appears in [24, p.338, 3.354.1] for $a, \mu > 0$

$$\int_0^\infty \frac{e^{-\mu v}}{a^2 + v^2} dv = \frac{1}{a}[\sin(a\mu)Ci(a\mu) - \cos(a\mu)si(a\mu)]$$

and this is also valid in the limit as $\mu \to 0$ because we easily see that

$$\int_0^\infty \frac{dv}{a^2 + v^2} = \frac{\pi}{2a} = -\frac{1}{a}si(0)$$

The substitution $v = -\log x$ in (5.7) gives us

$$\int_0^\infty \frac{e^{-\mu v}}{a^2 + v^2} dv = \int_0^1 \frac{t^{\mu-1}}{a^2 + \log^2 t} dt$$

and we therefore obtain

(5.10) $$\int_0^1 \frac{t^{\mu-1}}{a^2 + \log^2 t} dt = \frac{1}{a}[\sin(a\mu)Ci(a\mu) - \cos(a\mu)si(a\mu)]$$

With $\mu = 1$ we obtain the known result [24, p.523, 4.213.1]

$$(5.11) \quad \int_0^1 \frac{dt}{a^2 + \log^2 t} = \frac{1}{a}[\sin(a)Ci(a) - \cos(a)si(a)]$$

and differentiating this with respect to $a$ results in [24, p.523, 4.213.5].

Differentiation of (5.10) with respect to $\mu$ results in

$$a\int_0^1 \frac{t^{\mu-1} \log t}{a^2 + \log^2 t} dt = \sin(a\mu)\frac{\cos(a\mu)}{\mu} + a\cos(a\mu)Ci(a\mu) - \cos(a\mu)\frac{\sin(a\mu)}{\mu} + a\sin(a\mu)si(a\mu)$$

where we have used the elementary derivatives

$$\frac{d}{dx} si(tx) = \frac{\sin(tx)}{x} \quad \text{and} \quad \frac{d}{dx} Ci(tx) = \frac{\cos(tx)}{x}.$$

so that

$$(5.12) \quad \int_0^1 \frac{t^{\mu-1} \log t}{a^2 + \log^2 t} dt = \cos(a\mu)Ci(a\mu) + \sin(a\mu)si(a\mu)$$

With $\mu = 1$ we obtain the known result [24, p.523, 4.213.3]

$$(5.13) \quad \int_0^1 \frac{\log t}{a^2 + \log^2 t} dt = \cos(a)Ci(a) + \sin(a)si(a)$$

and differentiating this with respect to $a$ results in [24, p.524, 4.213.7]

$$(5.14) \quad \int_0^1 \frac{\log t}{\left[a^2 + \log^2 t\right]^2} dt = \frac{1}{2a}[\sin(a)Ci(a) - \cos(a)si(a)] - \frac{1}{2a^2}$$

A further differentiation of (5.12) with respect to $u$ gives us

$$(5.15) \quad \int_0^1 \frac{t^{\mu-1} \log^2 t}{a^2 + \log^2 t} dt = \frac{1}{\mu} - a[\sin(a\mu)Ci(a\mu) - \cos(a\mu)si(a\mu)]$$

which concurs with (5.10) because we note that

$$\int_0^1 \frac{t^{\mu-1}\log^2 t}{a^2 + \log^2 t}\,dt = \int_0^1 t^{\mu-1}\left[1 - \frac{a^2}{a^2 + \log^2 t}\right]dt$$

$\square$

Letting $a \to 2n\pi a$ in (5.10) gives us

$$\int_0^1 \frac{x^{\mu-1}dx}{4\pi^2 n^2 a^2 + \log^2 x} = \frac{1}{2n\pi a}[\sin(2n\pi a\mu)Ci(2n\pi a\mu) - \cos(2n\pi a\mu)si(2n\pi a\mu)]$$

and we make the summation

$$\sum_{n=1}^{\infty}\int_0^1 \frac{x^{\mu-1}dx}{4\pi^2 n^2 a^2 + \log^2 x} = \sum_{n=1}^{\infty}\frac{1}{2n\pi a}[\sin(2n\pi a\mu)Ci(2n\pi a\mu) - \cos(2n\pi a\mu)si(2n\pi a\mu)]$$

Assuming that interchanging the order of summation and integration is valid we have

$$\sum_{n=1}^{\infty}\int_0^1 \frac{x^{\mu-1}dx}{4\pi^2 n^2 a^2 + \log^2 x} = \frac{1}{a^2}\sum_{n=1}^{\infty}\int_0^1 \frac{x^{\mu-1}dx}{4\pi^2 n^2 + (\log^2 x)/a^2}$$

$$= \frac{1}{a^2}\int_0^1 \sum_{n=1}^{\infty}\frac{x^{\mu-1}dx}{4\pi^2 n^2 + (\log^2 x)/a^2}\,dx$$

We then use the well known identity ([11, p.296], [26, p.378])

(5.16) $$2\sum_{n=1}^{\infty}\frac{c}{4\pi^2 n^2 + c^2} = \frac{1}{e^c - 1} - \frac{1}{c} + \frac{1}{2}$$

to give us

$$2\sum_{n=1}^{\infty}\frac{1}{4\pi^2 n^2 + (\log^2 x)/a^2} = \left[\frac{1}{x^{1/a} - 1} - \frac{a}{\log x} + \frac{1}{2}\right]\frac{a}{\log x}$$

Hence we have

$$\frac{1}{a^2}\int_0^1 \sum_{n=1}^{\infty}\frac{x^{\mu-1}dx}{4\pi^2 n^2 + \log^2 x/a^2}\,dx = \frac{1}{2a}\int_0^1 \left[\frac{1}{x^{1/a} - 1} - \frac{a}{\log x} + \frac{1}{2}\right]\frac{x^{\mu-1}dx}{\log x}$$

$$= \sum_{n=1}^{\infty}\frac{1}{2n\pi a}[\sin(2n\pi a\mu)Ci(2n\pi a\mu) - \cos(2n\pi a\mu)si(2n\pi a\mu)]$$

and, at this stage, we recall (2.11)

$$\log \Gamma(a\mu) = \frac{1}{2}\log(2\pi) + \left(a\mu - \frac{1}{2}\right)\log(a\mu) - a\mu$$

$$+ \sum_{n=1}^{\infty} \frac{\sin(2n\pi a\mu)Ci(2n\pi a\mu) - \cos(2n\pi a\mu)si(2n\pi a\mu)}{n\pi}$$

so that we have

(5.17) $$\int_0^1 \left[\frac{1}{x^{1/a}-1} - \frac{a}{\log x} + \frac{1}{2}\right]\frac{x^{\mu-1}dx}{\log x} = \log \Gamma(a\mu) - \frac{1}{2}\log(2\pi) - \left(a\mu - \frac{1}{2}\right)\log(a\mu) + a\mu$$

With $a = 1$ we obtain [24, p.550, 4.283.10]

(5.18) $$\int_0^1 \left[\frac{1}{x-1} - \frac{1}{\log x} + \frac{1}{2}\right]\frac{x^{\mu-1}dx}{\log x} = \log \Gamma(\mu) - \frac{1}{2}\log(2\pi) - \left(\mu - \frac{1}{2}\right)\log \mu + \mu$$

With $\mu = 1$ in (5.18) we obtain [24, p.550, 4.283.2]

(5.19) $$\int_0^1 \left[\frac{1}{1-x} + \frac{1}{\log x} - \frac{1}{2}\right]\frac{dx}{\log x} = \frac{1}{2}\log(2\pi) - 1$$

With $\mu = 1$ in (5.17) we obtain

(5.20) $$\int_0^1 \left[\frac{1}{x^{1/a}-1} - \frac{a}{\log x} + \frac{1}{2}\right]\frac{dx}{\log x} = \log \Gamma(a) - \frac{1}{2}\log(2\pi) - \left(a - \frac{1}{2}\right)\log(a) + a$$

The substitution $u = x^{1/a}$ in (5.20) gives us

$$\int_0^1 \left[\frac{1}{x^{1/a}-1} - \frac{a}{\log x} + \frac{1}{2}\right]\frac{dx}{\log x} = \int_0^1 \left[\frac{1}{u-1} - \frac{1}{\log u} + \frac{1}{2}\right]\frac{u^{a-1}du}{\log u}$$

$$= \int_0^1 \left[\frac{u+1}{2(u-1)} - \frac{1}{\log u}\right]\frac{u^{a-1}du}{\log u}$$

Hence we obtain a slightly different derivation of [24, p.550, 4.283.10]

$$\int_0^1 \left[\frac{u+1}{2(u-1)} - \frac{1}{\log u}\right]\frac{u^{a-1}du}{\log u} = \log \Gamma(a) - \frac{1}{2}\log(2\pi) - \left(a - \frac{1}{2}\right)\log a + a$$

With $a = 1/2$ in (5.17) we obtain

$$(5.21) \quad \int_0^1 \left[\frac{1}{x^2-1} - \frac{1}{2\log x} + \frac{1}{2}\right] \frac{x^{\mu-1}dx}{\log x} = \log \Gamma\left(\frac{\mu}{2}\right) - \frac{1}{2}\log(2\pi) + \frac{1}{2}(\mu-1)\log(2\mu) + \frac{1}{2}\mu$$

and letting $\mu = 1$ in this results in [24, p.550, 4.283.5]

$$(5.22) \quad \int_0^1 \left[\frac{1}{x^2-1} - \frac{1}{2\log x} + \frac{1}{2}\right] \frac{dx}{\log x} = \frac{1}{2}[\log 2 - 1]$$

With the substitution $u = -\log x$ in (5.19) we get

$$\int_0^1 \left[\frac{1}{1-x} + \frac{1}{\log x} - \frac{1}{2}\right] \frac{dx}{\log x} = -\int_0^\infty \left[\frac{1}{e^u - 1} - \frac{1}{u} + \frac{1}{2}\right] \frac{e^{-u}}{u} du$$

and we thereby obtain the known integral [38, p.16]

$$(5.23) \quad \int_0^\infty \left[\frac{1}{e^u - 1} - \frac{1}{u} + \frac{1}{2}\right] \frac{e^{-u}}{u} du = 1 - \frac{1}{2}\log(2\pi)$$

Subtracting (5.22) from (5.18) gives us

$$(5.24) \quad \int_0^1 \left[\frac{1}{u-1} - \frac{1}{\log u} + \frac{1}{2}\right] \frac{(u^{a-1}-1)du}{\log u} = \log \Gamma(a) - \left(a - \frac{1}{2}\right)\log a + a - 1$$

We now refer to Malmstén's formula [38, p.16] (which is also derived in equation (E.22g) of [16])

$$(5.25) \quad \log \Gamma(a) = \int_0^\infty \left[(a-1)e^{-t} - \frac{e^{-t} - e^{-at}}{1 - e^{-t}}\right] \frac{dt}{t}$$

and a change of variables $u = e^{-t}$ gives us

$$(5.26) \quad \log \Gamma(a) = \int_0^1 \left[\frac{u^{a-1} - 1}{u - 1} - a + 1\right] \frac{du}{\log u}$$

Subtracting this from (5.24) gives us

$$\int_0^1 \left[ a-1-\frac{u^{a-1}-1}{\log u}+\frac{u^{a-1}-1}{2}\right]\frac{du}{\log u}=-\left(a-\frac{1}{2}\right)\log a+a-1$$

Differentiating this results in

(5.27) $$\int_0^1 \left[1-u^{a-1}+\frac{1}{2}u^{a-1}\log u\right]\frac{du}{\log u}=\frac{1}{2a}-\log a$$

and we see that

$$\int_0^1 \left[1-u^{a-1}+\frac{1}{2}u^{a-1}\log u\right]\frac{du}{\log u}=\int_0^1 \frac{1-u^{a-1}}{\log u}du+\frac{1}{2a}$$

The last integral is well known and is easily evaluated as follows.

We have $1/r=\int_0^\infty e^{-rx}dx$ ($r>0$) and, integrating that expression, we obtain Frullani's integral

$$\int_1^a \frac{dr}{r}=\int_1^a dr \int_0^\infty e^{-rx}dx = \int_0^\infty dx \int_1^a e^{-rx}dr$$

which implies that

$$\log a = \int_0^\infty \frac{e^{-x}-e^{-ax}}{x}dx$$

The substitution $u=e^{-x}$ gives us the well known integral

(5.28) $$\log a = \int_0^1 \frac{u^{a-1}-1}{\log u}du$$

and this confirms the above analysis.

$\square$

We have from (5.12)

$$\int_0^1 \frac{x^{\mu-1}\log x}{4\pi^2 n^2 a^2 + \log^2 x}dx = \cos(2n\pi a\mu)Ci(2n\pi a\mu)+\sin(2n\pi a\mu)si(2n\pi a\mu)$$

and we make the summation

$$\sum_{n=1}^{\infty}\int_0^1 \frac{x^{\mu-1}\log x}{4\pi^2 n^2 a^2+\log^2 x}\,dx = \sum_{n=1}^{\infty}[\cos(2n\pi a\mu)Ci(2n\pi a\mu)+\sin(2n\pi a\mu)si(2n\pi a\mu)]$$

$$=\frac{1}{2}\left[\psi(a\mu)-\log(a\mu)+\frac{1}{2a\mu}\right]$$

where we have used (4.3).

Assuming that interchanging the order of summation and integration is valid we have

$$\sum_{n=1}^{\infty}\int_0^1 \frac{x^{\mu-1}\log x}{4\pi^2 n^2 a^2+\log^2 x}\,dx = \frac{1}{a}\int_0^1 \sum_{n=1}^{\infty}\frac{x^{\mu-1}(\log x)/a}{4\pi^2 n^2+(\log^2 x)/a^2}\,dx$$

and using (5.16) gives us

$$\sum_{n=1}^{\infty}\frac{x^{\mu-1}(\log x)/a}{4\pi^2 n^2+(\log^2 x)/a^2} = \frac{1}{2}\left[\frac{1}{x^{1/a}-1}-\frac{a}{\log x}+\frac{1}{2}\right]x^{\mu-1}$$

so that

$$\psi(a\mu)-\log(a\mu)+\frac{1}{2a\mu} = \frac{1}{a}\int_0^1\left[\frac{1}{x^{1/a}-1}-\frac{a}{\log x}+\frac{1}{2}\right]x^{\mu-1}dx$$

or equivalently

(5.29) $$\psi(a\mu)-\log(a\mu) = \frac{1}{a}\int_0^1\left[\frac{1}{x^{1/a}-1}-\frac{a}{\log x}\right]x^{\mu-1}dx$$

The substitution $u = x^{1/a}$ results in

$$\frac{1}{a}\int_0^1\left[\frac{1}{x^{1/a}-1}-\frac{a}{\log x}\right]x^{\mu-1}dx = \int_0^1\left[\frac{1}{u-1}-\frac{1}{\log u}\right]u^{a\mu-1}du$$

so that

$$\psi(a\mu)-\log(a\mu) = \int_0^1\left[\frac{1}{u-1}-\frac{1}{\log u}\right]u^{a\mu-1}dx$$

Using (5.28) we obtain the well known integral representation

$$\psi(a) = \int_0^1 \left[ \frac{u^{a-1}}{u-1} - \frac{1}{\log u} \right] dx$$

With $a = 1$ we obtain the known result [9, p.178]

$$\gamma = \int_0^1 \left[ \frac{1}{1-u} + \frac{1}{\log u} \right] du$$

Direct integration shows that

$$\int_0^1 \left[ \frac{1}{1-u} + \frac{1}{\log u} \right] du = -\log(1-u) + li(u) \Big|_0^1$$

where $li(u)$ is the logarithm function defined by $li(u) = \int_0^u \frac{dt}{\log t}$. Hence we see that

(5.30) $\quad \gamma = \lim_{u \to 1} [li(u) - \log(1-u)]$

$\square$

We also have from (5.7)

$$\int_0^\infty \frac{e^{-\mu v}}{n^2 + v^2} dv = \frac{1}{n} [\sin(n\mu) Ci(n\mu) - \cos(n\mu) si(n\mu)]$$

and we make the summation

$$\sum_{n=1}^\infty \int_0^\infty \frac{e^{-\mu v}}{n^2 + v^2} dv = \sum_{n=1}^\infty \frac{1}{n} [\sin(n\mu) Ci(n\mu) - \cos(n\mu) si(n\mu)]$$

Assuming that interchanging the order of summation and integration is valid

$$\sum_{n=1}^\infty \int_0^\infty \frac{e^{-\mu v}}{n^2 + v^2} dv = \int_0^\infty e^{-\mu v} \sum_{n=1}^\infty \frac{1}{n^2 + v^2} dv$$

and using

(5.31) $\quad \pi \coth \pi v = \frac{1}{v} + 2v \sum_{n=1}^\infty \frac{1}{n^2 + v^2}$

we have

(5.32) $$\frac{1}{2}\int_0^\infty \frac{\pi v \coth \pi v - 1}{v^2} e^{-\mu v} dv = \sum_{n=1}^\infty \frac{1}{n}[\sin(n\mu)Ci(n\mu) - \cos(n\mu)si(n\mu)]$$

We note that *Mathematica* cannot evaluate this integral.

With $\mu = \pi$ we have

$$\frac{1}{2}\int_0^\infty \frac{\pi v \coth \pi v - 1}{v^2} e^{-\pi v} dv = \sum_{n=1}^\infty \frac{(-1)^n}{n} si(n\pi)$$

and since

$$\sum_{n=1}^\infty \frac{(-1)^n}{n} si(n\pi) = \frac{\pi}{2}\log 2 - \frac{\pi}{2}$$

we obtain

(5.33) $$\frac{1}{2}\int_0^\infty \frac{\pi v \coth \pi v - 1}{v^2} e^{-\pi v} dv = \frac{\pi}{2}\log 2 - \frac{\pi}{2}$$

With $\mu = 2\pi$ we have

$$\frac{1}{2}\int_0^\infty \frac{\pi v \coth \pi v - 1}{v^2} e^{-2\pi v} dv = \sum_{n=1}^\infty \frac{si(2n\pi)}{n}$$

and since

$$\sum_{n=1}^\infty \frac{si(2n\pi)}{n} = \frac{\pi}{2}\log(2\pi) - \pi$$

we obtain

(5.34) $$\frac{1}{2}\int_0^\infty \frac{\pi v \coth \pi v - 1}{v^2} e^{-2\pi v} dv = \frac{\pi}{2}\log(2\pi) - \pi$$

More generally, using (2.11), we see that

(5.34.1) $$\frac{1}{2\pi}\int_0^\infty \frac{\pi v \coth \pi v - 1}{v^2} e^{-\mu v} dv = \log \Gamma(\mu) - \frac{1}{2}\log(2\pi) - \left(\mu - \frac{1}{2}\right)\log \mu + \mu$$

□

The substitution $u = \tan x$ gives us

$$\int_0^{\pi/2} \exp(-p \tan x)\,dx = \int_0^\infty \frac{e^{-pu}}{1+u^2}\,du$$

$$= \sin p\, Ci(p) - \cos p\, si(p)$$

where we have used (5.7). We therefore obtain [24, p.336, 3.341]

(5.35) $$\int_0^{\pi/2} \exp(-p \tan x)\,dx = \sin p\, Ci(p) - \cos p\, si(p)$$

so that

$$\int_0^{\pi/2} \exp(-2n\pi a \tan x)\,dx = \sin 2n\pi a\, Ci(2n\pi a) - \cos 2n\pi a\, si(2n\pi a)$$

Referring to (2.11)

$$\log \Gamma(a) = \frac{1}{2}\log(2\pi) + \left(a - \frac{1}{2}\right)\log a - a + \frac{1}{\pi}\sum_{n=1}^\infty \frac{1}{n}[\sin(2n\pi a)Ci(2n\pi a) - \cos(2n\pi a)si(2n\pi a)]$$

and using

$$\sum_{n=1}^\infty \frac{1}{n}\exp(-2n\pi a \tan x) = -\log[1 - \exp(-2\pi a \tan x)]$$

we obtain an integral representation for $\log \Gamma(a)$

(5.36) $$\log \Gamma(a) = \frac{1}{2}\log(2\pi) + \left(a - \frac{1}{2}\right)\log a - a - \frac{1}{\pi}\int_0^{\pi/2} \log[1 - \exp(-2\pi a \tan x)]\,dx$$

These integrals complement the related ones appearing in [29].

**6. Binet's integrals for the digamma function**

We substitute $\mu = n$ in (5.9) and make the summation

$$\sum_{n=1}^\infty \int_0^\infty \frac{ve^{-nv}}{a^2 + v^2}\,dv = -\sum_{n=1}^\infty [\cos(na)Ci(na) + \sin(na)si(na)]$$

and the geometric series gives us

$$\sum_{n=1}^{\infty}\int_0^{\infty}\frac{ve^{-nv}}{a^2+v^2}\,dv = \int_0^{\infty}\frac{ve^{-v}}{(a^2+v^2)(1-e^{-v})}\,dv$$

so that we have

$$\int_0^{\infty}\frac{v}{(a^2+v^2)(e^v-1)}\,dv = -\sum_{n=1}^{\infty}[Ci(na)\cos(na)+si(na)\sin(na)]$$

We now let $a \to 2\pi a$ and $v \to 2\pi t$ to obtain

(6.1) $$\int_0^{\infty}\frac{t}{(a^2+t^2)(e^{2\pi t}-1)}\,dt = -\sum_{n=1}^{\infty}[\cos(2n\pi a)Ci(2n\pi a)+\sin(2n\pi a)si(2n\pi a)]$$

and using (4.3)

$$\psi(a) = \log a - \frac{1}{2a} + 2\sum_{n=1}^{\infty}[\cos(2n\pi a)Ci(2n\pi a)+\sin(2n\pi a)si(2n\pi a)]$$

we then see that

(6.2) $$\psi(a) = \log a - \frac{1}{2a} - 2\int_0^{\infty}\frac{t}{(a^2+t^2)(e^{2\pi t}-1)}\,dt$$

Equation (6.2) is well known and, inter alia, is reported in [40, p.251] and is frequently obtained by differentiating Binet's second formula for the log gamma function.

Alternatively, we substitute $a = 2\pi n$ in (5.9) and make the summation

$$\sum_{n=1}^{\infty}\int_0^{\infty}\frac{ve^{-\mu v}}{4\pi^2 n^2+v^2}\,dv = -\sum_{n=1}^{\infty}[\cos(2n\pi\mu)Ci(2n\pi\mu)+\sin(2n\pi\mu)si(2n\pi\mu)]$$

Assuming that it is valid to change the order of summation and integration

$$\sum_{n=1}^{\infty}\int_0^{\infty}\frac{ve^{-\mu v}}{4\pi^2 n^2+v^2}\,dv = \int_0^{\infty}\sum_{n=1}^{\infty}\frac{ve^{-\mu v}}{4\pi^2 n^2+v^2}\,dv$$

and using (5.16) we have

$$\sum_{n=1}^{\infty}\int_0^{\infty}\frac{ve^{-\mu v}}{4\pi^2 n^2+v^2}dv=\frac{1}{2}\int_0^{\infty}\left[\frac{1}{e^v-1}-\frac{1}{v}+\frac{1}{2}\right]e^{-\mu v}dv$$

Therefore we obtain

$$\frac{1}{2}\int_0^{\infty}\left[\frac{1}{e^v-1}-\frac{1}{v}+\frac{1}{2}\right]e^{-\mu v}dv=-\sum_{n=1}^{\infty}[\cos(2n\pi\mu)Ci(2n\pi\mu)+\sin(2n\pi\mu)si(2n\pi\mu)]$$

and using (4.3) this results in

(6.3) $$\int_0^{\infty}\left[\frac{1}{e^v-1}-\frac{1}{v}+\frac{1}{2}\right]e^{-\mu v}dv=\log\mu-\psi(\mu)-\frac{1}{2\mu}$$

which may be obtained by differentiating Binet's first formula for the log gamma function.

## 7. The Barnes double gamma function

Integration by parts results in

$$\int x[\cos(ax)Ci(ax)+\sin(ax)si(ax)]dx=\frac{x}{a}[\sin(ax)Ci(ax)-\cos(ax)si(ax)]$$

$$-\frac{1}{a}\int[\sin(ax)Ci(ax)-\cos(ax)si(ax)]dx$$

where we note that

$$\frac{d}{dx}\frac{1}{a}[\sin(ax)Ci(ax)-\cos(ax)si(ax)]=\cos(ax)Ci(ax)+\sin(ax)si(ax)$$

because $\frac{d}{dx}si(ax)=\frac{\sin(ax)}{x}$ and $\frac{d}{dx}Ci(ax)=\frac{\cos(ax)}{x}$.

We have

$$\int[\sin(ax)Ci(ax)-\cos(ax)si(ax)]dx=-\frac{1}{a}[\cos(ax)Ci(ax)+\sin(ax)si(ax)-\log x]$$

which may be verified by noting that

$$\frac{d}{dx}\frac{1}{a}[\cos(ax)Ci(ax)+\sin(ax)si(ax)-\log x]=\sin(ax)Ci(ax)-\cos(ax)si(ax)$$

Hence we obtain

$$\int x[\cos(ax)Ci(ax)+\sin(ax)si(ax)]dx = \frac{x}{a}[\sin(ax)Ci(ax)-\cos(ax)si(ax)]$$

$$+\frac{1}{a^2}[\cos(ax)Ci(ax)+\sin(ax)si(ax)-\log x]$$

and we have the definite integral

(7.1) $$\int_0^u x[\cos(ax)Ci(ax)+\sin(ax)si(ax)]dx = \frac{u}{a}[\sin(au)Ci(au)-\cos(au)si(au)]$$

$$+\frac{1}{a^2}[\cos(au)Ci(au)+\sin(au)si(au)-\log u - \gamma - \log a]$$

where we have employed the limits (2.5) and (2.6)

$$\lim_{u \to 0}[\cos(au)Ci(au) - \log u] = \gamma + \log a$$

$$\lim_{u \to 0}\sin(au)Ci(au) = 0$$

From (4.3) we see that

$$x\psi(x) - x\log x + \frac{1}{2} = 2\sum_{n=1}^{\infty} x[\cos(2n\pi x)Ci(2n\pi x)+\sin(2n\pi x)si(2n\pi x)]$$

and we have the integral

$$\frac{1}{2}\int_0^u \left[x\psi(x) - x\log x + \frac{1}{2}\right]dx = \sum_{n=1}^{\infty}\int_0^u x[\cos(2n\pi x)Ci(2n\pi x)+\sin(2n\pi x)si(2n\pi x)]dx$$

$$= u\sum_{n=1}^{\infty} \frac{\sin(2n\pi u)Ci(2n\pi u) - \cos(2n\pi u)si(2n\pi u)}{2n\pi}$$

$$+\sum_{n=1}^{\infty} \frac{\cos(2n\pi u)Ci(2n\pi u)+\sin(2n\pi u)si(2n\pi u) - \log u - \gamma - \log(2n\pi)}{(2n\pi)^2}$$

and, using (2.11) this may be written as

$$= \frac{u}{2}\left[\log\Gamma(u) - \frac{1}{2}\log(2\pi) - \left(u - \frac{1}{2}\right)\log u + u\right]$$

$$+ \sum_{n=1}^{\infty}\frac{\cos(2n\pi u)\,Ci(2n\pi u) + \sin(2n\pi u)\,si(2n\pi u) - \log u - \gamma - \log(2n\pi)}{(2n\pi)^2}$$

$$= \frac{u}{2}\left[\log\Gamma(u) - \frac{1}{2}\log(2\pi) - \left(u - \frac{1}{2}\right)\log u + u\right]$$

$$+ \sum_{n=1}^{\infty}\frac{\cos(2n\pi u)\,Ci(2n\pi u) + \sin(2n\pi u)\,si(2n\pi u)}{(2n\pi)^2}$$

$$- \frac{[\gamma + \log(2\pi u)]\varsigma(2) - \varsigma'(2)}{4\pi^2}$$

We note that

$$\int_0^u \left[x\psi(x) - x\log x + \frac{1}{2}\right]dx = \int_0^u x\psi(x)\,dx - \frac{u^2}{4}[2\log u - 1] + \frac{u}{2}$$

Integration by parts results in

$$\int_0^u x\psi(x)\,dx = u\log\Gamma(u) - \int_0^u \log\Gamma(x)\,dx$$

where we have used $\lim_{u\to 0}[u\log\Gamma(u)] = \lim_{u\to 0}[u\log\Gamma(1+u) - u\log u] = 0$.

Using Alexeiewsky's theorem [38, p.32]

(7.2) $\quad \int_0^u \log\Gamma(x)\,dx = u\log\Gamma(u) - \log G(1+u) - \frac{1}{2}u(u-1) + \frac{1}{2}u\log(2\pi)$

where the Barnes double gamma function $\Gamma_2(x) = 1/G(x)$ is defined by [38, p.25]

(7.3) $\quad G(1+u) = (2\pi)^{u/2}\exp\left[-\frac{1}{2}(\gamma u^2 + u^2 + u)\right]\prod_{k=1}^{\infty}\left\{\left(1 + \frac{u}{k}\right)^k \exp\left(\frac{u^2}{2k} - u\right)\right\}$

we obtain

$$\int_0^u x\psi(x)dx = \log G(1+u) + \frac{1}{2}u^2 - \frac{1}{2}u - \frac{1}{2}u\log(2\pi)$$

Therefore we have

$$\int_0^u \left[x\psi(x) - x\log x + \frac{1}{2}\right]dx = \log G(1+u) - \frac{1}{2}u\log(2\pi) - \frac{u^2}{2}\left[\log u - \frac{3}{2}\right]$$

and hence we see that

$$\sum_{n=1}^\infty \frac{\cos(2n\pi u)\,Ci(2n\pi u) + \sin(2n\pi u)\,si(2n\pi u)}{(2n\pi)^2}$$

$$= \frac{1}{2}[\log G(1+u) - u\log\Gamma(u)] + \frac{1}{4}u(u-1)\log u - \frac{1}{8}u^2 + \frac{[\gamma + \log(2\pi u)]\varsigma(2) - \varsigma'(2)}{4\pi^2}$$

Using the functional equation for the Riemann zeta function it is easily shown that

$$\varsigma'(-1) = \frac{1}{12}(1 - \gamma - \log 2\pi) + \frac{1}{2\pi^2}\varsigma'(2)$$

which may be written as

(7.4) $$\quad \frac{1}{2}\left[\frac{1}{12} - \varsigma'(-1)\right] = \frac{[\gamma + \log(2\pi)]\varsigma(2) - \varsigma'(2)}{4\pi^2}$$

Hence we obtain

(7.5) $$\quad \sum_{n=1}^\infty \frac{\cos(2n\pi u)\,Ci(2n\pi u) + \sin(2n\pi u)\,si(2n\pi u)}{(2n\pi)^2}$$

$$= \frac{1}{2}[\log G(1+u) - u\log\Gamma(u)] + \frac{1}{4}\left[u(u-1) + \frac{1}{6}\right]\log u - \frac{1}{8}u^2 + \frac{1}{2}\left[\frac{1}{12} - \varsigma'(-1)\right]$$

Letting $u = 0$ in (7.5) results in (after using the limit (2.5))

$$\sum_{n=1}^\infty \frac{\gamma + \log(2n\pi)}{(2n\pi)^2} = \frac{1}{2}\left[\frac{1}{12} - \varsigma'(-1)\right]$$

which is simply another representation of (7.4).

Letting $u = 1$ results in

$$(7.5) \quad \sum_{n=1}^{\infty} \frac{Ci(2n\pi)}{(2n\pi)^2} = -\frac{1}{12} - \frac{1}{2}\varsigma'(-1)$$

□

We note the Gosper/Vardi identity

$$(7.6) \quad \log G(u+1) - u \log \Gamma(u) = \varsigma'(-1) - \varsigma'(-1, u)$$

This functional equation was derived by Vardi in 1988 and also by Gosper in 1997 (see [2]). Then, using the well-known result [5, p.264] in terms of the Bernoulli polynomials

$$\varsigma(1-n, u) = -\frac{B_n(u)}{n} \quad \text{for } n \geq 1$$

we obtain

$$(7.7) \quad \varsigma'(-1, u) =$$

$$-\varsigma(-1, u) \log u - \frac{1}{4}u^2 + \frac{1}{12} - \frac{1}{2\pi^2} \sum_{n=1}^{\infty} \frac{1}{n^2} [\cos(2n\pi u) Ci(2n\pi u) + \sin(2n\pi u) si(2n\pi u)]$$

which was reported by Elizalde [21] in 1985 for $u > 0$.

□

Integrating (6.1) gives us

$$(7.8) \quad \int_0^x da \int_0^{\infty} \frac{avdv}{(a^2+v^2)(e^{2\pi v}-1)} = -\sum_{n=1}^{\infty} \int_0^x a[\cos(2n\pi a) Ci(2n\pi a) + \sin(2n\pi a) si(2n\pi a)] da$$

and using (7.1) we obtain

$$\sum_{n=1}^{\infty} \int_0^x a[\cos(2n\pi a) Ci(2n\pi a) + \sin(2n\pi a) si(2n\pi a)] da$$

$$= \sum_{n=1}^{\infty} \frac{\cos(2n\pi x) Ci(2n\pi x) + \sin(2n\pi x) si(2n\pi x)}{(2n\pi)^2} - \frac{1}{24}(\log x + \log(2\pi) + \gamma) + \frac{\varsigma'(2)}{4\pi^2}$$

$$+ x \sum_{n=1}^{\infty} \frac{\sin(2n\pi x) Ci(2n\pi x) - \cos(2n\pi x) si(2n\pi x)}{2n\pi}$$

Using (2.11) this becomes

$$= \sum_{n=1}^{\infty} \frac{\cos(2n\pi x)Ci(2n\pi x)+\sin(2n\pi x)si(2n\pi x)}{(2n\pi)^2} - \frac{1}{24}(\log x + \log(2\pi) + \gamma) + \frac{\varsigma'(2)}{4\pi^2}$$

$$+ \frac{1}{2}x\log\Gamma(x) - \frac{x}{4}\log(2\pi) - \frac{x}{2}\left(x - \frac{1}{2}\right)\log x + \frac{1}{2}x^2$$

Alternatively, integrating the left-hand side of (7.8) gives us

$$I = \int_0^x da \int_0^\infty \frac{av}{(a^2+v^2)(e^{2\pi v}-1)} dv = \frac{1}{2}\int_0^\infty \frac{v\, dv}{e^{2\pi v}-1}\int_0^x \frac{2a}{a^2+v^2} da$$

$$= \frac{1}{2}\int_0^\infty \frac{v\log(x^2+v^2)}{e^{2\pi v}-1} dv - \int_0^\infty \frac{v\log v}{e^{2\pi v}-1} dv$$

We have [17]

(7.9) $$\log G(1+x) = \frac{1}{2}x^2\left[\log x - \frac{3}{2}\right] + \frac{1}{2}x\log(2\pi) + \varsigma'(-1) - \int_0^\infty \frac{v\log(v^2+x^2)}{e^{2\pi v}-1} dv$$

where $G(x)$ is the Barnes double gamma function. This result was originally obtained by Adamchik [3] in 2004. With $x = 0$ we have

(7.10) $$\int_0^\infty \frac{v\log v}{e^{2\pi v}-1} dv = \frac{1}{2}\varsigma'(-1)$$

Therefore we obtain

$$I = \frac{1}{4}x^2\left[\log x - \frac{3}{2}\right] + \frac{1}{4}x\log(2\pi) - \frac{1}{2}\log G(1+x)$$

Hence we obtain

(7.11) $$\sum_{n=1}^{\infty} \frac{\cos(2n\pi x)Ci(2n\pi x)+\sin(2n\pi x)si(2n\pi x)}{(2n\pi)^2}$$

$$= \frac{1}{24}(\log x + \log(2\pi) + \gamma) - \frac{1}{4}x^2\left[\log x - \frac{3}{2}\right] - \frac{1}{2}x\log(2\pi) + \frac{1}{2}\log G(1+x) + \frac{\varsigma'(2)}{4\pi^2}$$

$$-\frac{1}{2}x\log\Gamma(x) + \frac{x}{2}\left(x - \frac{1}{2}\right)\log x - \frac{1}{2}x^2$$

which is equivalent to (7.5).

Alternatively, we multiply (6.2) by $a$ and integrate to obtain

$$\int_0^x a\psi(a)da = \frac{1}{2}x^2\log x - \frac{1}{4}x^2 - \frac{1}{2}x - \frac{1}{2}x^2\left[\log x - \frac{3}{2}\right] - \frac{1}{2}x\log(2\pi) + \log G(1+x)$$

Integration by parts results in

$$\int_0^x a\psi(a)da = x\log\Gamma(x) - \int_0^x \log\Gamma(a)da$$

and we obtain Alexeiewsky's theorem [38, p.32]

$$(7.12) \quad \int_0^x \log\Gamma(a)da = x\log\Gamma(x) - \log G(1+x) - \frac{1}{2}x^2 + \frac{1}{2}x + \frac{1}{2}x\log(2\pi)$$

$\square$

We may also derive (7.5) by employing (2.1) and making the summation

$$\sum_{n=1}^{\infty}\frac{1}{n}\int_0^1 \log\Gamma(x+a)\sin 2n\pi x\, dx = -\sum_{n=1}^{\infty}\frac{1}{2n^2\pi}\left[\log a - \cos(2n\pi a)Ci(2n\pi a) - \sin(2n\pi a)si(2n\pi a)\right]$$

Assuming that $\sum_{n=1}^{\infty}\frac{1}{n}\int_0^1 \log\Gamma(x+a)\sin 2n\pi x\, dx = \int_0^1 \log\Gamma(x+a)\sum_{n=1}^{\infty}\frac{\sin 2n\pi x}{n}dx$ we have

$$\int_0^1 \log\Gamma(x+a)\sum_{n=1}^{\infty}\frac{\sin 2n\pi x}{n}dx = -\sum_{n=1}^{\infty}\frac{1}{2n^2\pi}\left[\log a - \cos(2n\pi a)Ci(2n\pi a) - \sin(2n\pi a)si(2n\pi a)\right]$$

and substituting $\sum_{n=1}^{\infty}\frac{\sin 2n\pi x}{n} = \frac{\pi}{2}(1-2x)$ gives us

(7.13)

$$\int_0^1 (1-2x)\log\Gamma(x+a)\,dx = -\frac{1}{6}\log a + \frac{1}{\pi^2}\sum_{n=1}^{\infty}\frac{\cos(2n\pi a)Ci(2n\pi a)+\sin(2n\pi a)si(2n\pi a)}{n^2}$$

Barnes [38, p.207] determined that

$$\int_0^z \log\Gamma(x+a)\,dx = \left(\frac{1}{2}+\frac{1}{2}\log(2\pi)-a\right)z - \frac{1}{2}z^2$$

$$+(z+a-1)\log\Gamma(z+a)-\log G(z+a)+(1-a)\log\Gamma(a)+\log G(a)$$

so that

(7.14) $\quad \int_0^1 \log\Gamma(x+a)\,dx = \frac{1}{2}\log(2\pi)+a\log a - a$

which we derived in a very different manner at the start of this paper.

We have from [38, p.209]

$$2\int_0^z x\log\Gamma(x+a)\,dx = \left(\frac{1}{4}-\frac{1}{2}a+\frac{1}{2}a^2-2\log A\right)z + \left(\frac{1}{2}\log(2\pi)-\frac{1}{2}a+\frac{1}{4}\right)z^2$$

$$-\frac{1}{2}z^3 + z^2\log\Gamma(z+a)-(a-1)^2\left[\log\Gamma(z+a)-\log\Gamma(a)\right]+(2a-3)\left[\log G(z+a)-\log G(a)\right]$$

$$-2\left[\log\Gamma_3(z+a)-\log\Gamma_3(a)\right]$$

so that

$$2\int_0^1 t\log\Gamma(a+t)\,dt = \left(\frac{1}{4}-\frac{1}{2}a+\frac{1}{2}a^2-2\log A\right) + \left(\frac{1}{2}\log(2\pi)-\frac{1}{2}a+\frac{1}{4}\right)$$

$$-\frac{1}{2}+\log\Gamma(1+a)-(a-1)^2\left[\log\Gamma(1+a)-\log\Gamma(a)\right]+(2a-3)\left[\log G(1+a)-\log G(a)\right]$$

$$-2\left[\log\Gamma_3(1+a)-\log\Gamma_3(a)\right]$$

We have

$$\log A = \frac{1}{12} - \varsigma'(-1)$$

where $A$ is the Glaisher-Kinkelin constant and from [38, p.25] the numerical value is reported as approximately

$$A = 1.282427130...$$

Since $\log \Gamma_3(1+a) = \log \Gamma_2(a) + \log \Gamma_3(a) = -\log G(a) + \log \Gamma_3(a)$ we obtain

$$(7.15)\, 2\int_0^1 t \log \Gamma(a+t)\, dt = -a + \frac{1}{2}a^2 - 2\log A + \frac{1}{2}\log(2\pi)$$

$$+ \log a - (a-1)^2 \log a + 2(a-1)\log \Gamma(a) - 2\log G(a)$$

Therefore, using (7.13), (7.14) and (7.15), we again obtain (7.5).

□

We see from (5.9) that

$$\sum_{n=1}^{\infty} \frac{1}{n^2} \int_0^{\infty} \frac{v e^{-nv}}{a^2 + v^2}\, dv = -\sum_{n=1}^{\infty} \frac{1}{n^2}[\cos(na)Ci(na) + \sin(na)si(na)]$$

and we have

$$\sum_{n=1}^{\infty} \frac{1}{n^2} \int_0^{\infty} \frac{v e^{-nv}}{(2\pi u)^2 + v^2}\, dv = \int_0^{\infty} \frac{v Li_2(e^{-v})}{(2\pi u)^2 + v^2}\, dv$$

where $Li_s(x)$ is the polylogarithm function $Li_s(x) = \sum_{n=1}^{\infty} \frac{x^n}{n^s}$.

We then obtain

$$-\frac{1}{4\pi^2} \int_0^{\infty} \frac{v Li_2(e^{-v})}{(2\pi u)^2 + v^2}\, dv = \frac{1}{2}[\log G(1+u) - u \log \Gamma(u)] + \frac{1}{4}\left[u(u-1) + \frac{1}{6}\right]\log u - \frac{1}{8}u^2$$

$$+ \frac{1}{2}\left[\frac{1}{12} - \varsigma'(-1)\right]$$

With the substitution $x = e^{-v}$ we have

$$\frac{1}{4\pi^2} \int_0^1 \frac{Li_2(x)}{(2\pi u)^2 + \log^2 x} \frac{\log x}{x}\, dv = \frac{1}{2}[\log G(1+u) - u \log \Gamma(u)] + \frac{1}{4}\left[u(u-1) + \frac{1}{6}\right]\log u - \frac{1}{8}u^2$$

$$+\frac{1}{2}\left[\frac{1}{12}-\varsigma'(-1)\right]$$

□

Further integrals and series involving the sine and cosine integrals are contained, inter alia, in [15], [18] and [19].

**8. Oscillatory integrals**

Let us assume that (i) the oscillatory integral $\int_0^\infty h(x)\sin(tx)\,dx$ exists where $t > 0$ and (ii) $h(x)$ is strictly decreasing for $x \in (0,\infty)$.

Then, as shown by Tuck [39], the following integral is strictly positive

(8.1) $$\int_0^\infty h(x)\sin(tx)\,dx > 0$$

Using (5.1) it is easily shown that

$$\int_0^\infty \frac{\sin tx}{a+x}\,dx = \sin(ta)Ci(ta) - \cos(ta)si(ta)$$

We note that $h(x) = \dfrac{1}{a+x}$ satisfies the conditions of (8.1) and we then deduce that

$$\sin(ta)Ci(ta) - \cos(ta)si(ta) > 0$$

This may also be deduced from (5.10).

It is shown in [18] that

$$\log\Gamma(a) = \frac{1}{2}\log(2\pi) + \left(a - \frac{1}{2}\right)\log a - a + \frac{1}{\pi}\sum_{n=1}^\infty \frac{1}{n}[\sin(2n\pi a)Ci(2n\pi a) - \cos(2n\pi a)si(2n\pi a)]$$

and we then deduce that

(8.2) $$\log\Gamma(a) > \frac{1}{2}\log(2\pi) + \left(a - \frac{1}{2}\right)\log a - a$$

□

Let us assume that (i) the integral $\int_0^\infty h(x)\cos(tx)\,dx$ exists where $t > 0$, (ii) $h(x) > 0$ for $x \in (0,\infty)$ and (iii) $h(x)$ is strictly decreasing for $x \in (0,\infty)$, (iv) $\lim_{x\to\infty} h(x) = 0$ and $\lim_{x\to 0+} x\,h(x) = 0$, (v) $h(x)$ is strictly convex and (v) the derivative $h'(x)$ exists and $\lim_{x\to\infty} h'(x) = 0$.

Then the following integral is strictly positive

$$(8.3) \qquad \int_0^\infty h(x)\cos(tx)\,dx > 0$$

which was also shown by Tuck [39].

We have

$$(8.4) \qquad \int_0^\infty \frac{\cos 2n\pi u}{1+u}\,du = -Ci(2n\pi)$$

and this enables us to deduce that $Ci(2n\pi) < 0$.

**9. Open access to our own work**

This paper contains references to a number of other papers and most of them are currently freely available on the internet. Surely now is the time that all of <u>our</u> work should be freely accessible by <u>all</u>. The mathematics community should lead the way on this by publishing <u>everything</u> on arXiv, or in an equivalent open access repository. We think it, we write it, so why hide it? You know it makes sense.

**REFERENCES**

[1]   M. Abramowitz and I.A. Stegun (Eds.), Handbook of Mathematical Functions with Formulas, Graphs and Mathematical Tables. Dover, New York, 1970.
http://www.math.sfu.ca/~cbm/aands/

[2]   V.S. Adamchik, Contributions to the Theory of the Barnes Function. Computer Physics Communications, 2003.
http://www-2.cs.cmu.edu/~adamchik/articles/barnes1.pdf

[3]   V.S. Adamchik, Symbolic and numeric computations of the Barnes function. Computer Physics Communications, 157 (2004) 181-190.
Symbolic and Numeric Computations of the Barnes Function

[4]   T.M. Apostol, Mathematical Analysis, Second Ed., Addison-Wesley Publishing

Company, Menlo Park (California), London and Don Mills (Ontario), 1974.

[5]  T.M. Apostol, Introduction to Analytic Number Theory.
     Springer-Verlag, New York, Heidelberg and Berlin, 1976.

[6]  R.G. Bartle, The Elements of Real Analysis.
     2nd Ed. John Wiley & Sons Inc., 1976.

[7]  G. Bateman and A. Erdélyi, Tables of Integral Transforms I.
     McGraw Hill Book Company, New York, 1954.

[8]  B.C. Berndt, The Gamma Function and the Hurwitz Zeta Function.
     Amer. Math. Monthly, 92,126-130, 1985.
     http://202.120.1.3/new/course/fourier/supplement/Zeta-Function.pdf

[9]  G. Boros and V.H. Moll, Irresistible Integrals: Symbolics, Analysis and
     Experiments in the Evaluation of Integrals. Cambridge University Press, 2004.

[10] D.M. Bradley, A class of series acceleration formulae for Catalan's constant.
     The Ramanujan Journal, Vol. 3, Issue 2, 159-173, 1999.
     http://germain.umemat.maine.edu/faculty/bradley/papers/rj.pdf

[11] T.J.I'a Bromwich, Introduction to the theory of infinite series. Third edition.
     AMS Chelsea Publishing, 1991.

[12] H.S. Carslaw, Introduction to the theory of Fourier Series and Integrals.
     Third Ed. Dover Publications Inc, 1930.

[13] M.W. Coffey, On one-dimensional digamma and polygamma series related to
     the evaluation of Feynman diagrams.
     J. Comput. Appl. Math, 183, 84-100, 2005. arXiv:math-ph/0505051 [pdf, ps, other]

[14] D.F. Connon, Some series and integrals involving the Riemann zeta function,
     binomial coefficients and the harmonic numbers. Volume II(b), 2007.
     arXiv:0710.4024 [pdf]

[15] D.F. Connon, Some series and integrals involving the Riemann zeta function,
     binomial coefficients and the harmonic numbers. Volume V, 2007.
     arXiv:0710.4047 [pdf]

[16] D.F. Connon, Some series and integrals involving the Riemann zeta function,
     binomial coefficients and the harmonic numbers. Volume VI, 2007.
     arXiv:0710.4032 [pdf]

[17] D.F. Connon, Some applications of the Stieltjes constants.
     arXiv:0901.2083 [pdf], 2009

[18]  D.F. Connon, Some trigonometric integrals involving the log gamma and the digamma function. 2010. arXiv:1005.3469 [pdf]

[19]  D.F. Connon, Some applications of the Dirichlet integrals to the summation of series and the evaluation of integrals involving the Riemann zeta function. (in preparation, 2012).

[20]  P.J. de Doelder, On some series containing $\psi(x)-\psi(y)$ and $(\psi(x)-\psi(y))^2$ for certain values of $x$ and $y$. J. Comput. Appl. Math. 37, 125-141, 1991.

[21]  E. Elizalde, Derivative of the generalised Riemann zeta function $\varsigma(z,q)$ at $z=-1$. J. Phys. A Math. Gen. (1985) 1637-1640

[22]  O. Espinosa and V.H. Moll, On some integrals involving the Hurwitz zeta function: Part I. The Ramanujan Journal, 6,150-188, 2002.
http://arxiv.org/abs/math.CA/0012078

[23] J.W.L. Glaisher, Tables of the Numerical Values of the Sine-Integral, Cosine-Integral, and Exponential-Integral.
Philosophical Transactions of the Royal Society of London
Vol. 160, (1870), pp. 367-388.

[24]  I.S. Gradshteyn and I.M. Ryzhik, Tables of Integrals, Series and Products.
Sixth Ed., Academic Press, 2000.
Errata for Sixth Edition http://www.mathtable.com/errata/gr6_errata.pdf

[25]  T.H. Gronwall, The gamma function in integral calculus.
Annals of Math., 20, 35-124, 1918.

[26]  K. Knopp, Theory and Application of Infinite Series.
Second English Ed. Dover Publications Inc, New York, 1990.

[27]  E.E. Kummer, Beitrag zur Theorie der Function $\Gamma(x) = \int_0^\infty e^{-v} v^{x-1} dv$.
J. Reine Angew. Math., 35, 1-4, 1847.
http://www.digizeitschriften.de/

[28]  M. Lerch, Sur la différentiation d'une classe de séries trigonométriques.
Annales scientifiques de l'École Normale Supérieure, Sér. 3,tome 12(1895), p.351-361.
http://www.numdam.org/

[29]  L.A. Medina and V.H. Moll, The integrals in Gradshteyn and Ryzhik.
Part 10: The digamma function. Scientia, Series A: Math. Sciences 17, 2009, 45-66.
http://129.81.170.14/~vhm/papers_html/final10.pdf

[30] M. Merkle and M. M. R. Merkle, Krull's theory for the double gamma function. Appl. Math. Comput., 218 (2011), 935-943, doi:10.1016/j.amc.2011.01.090. http://www.milanmerkle.com/documents/radovi/KRULLcs.pdf

[31] N. Nielsen, Theorie des Integrallogarithmus und verwanter tranzendenten. 1906. http://www.math.uni-bielefeld.de/~rehmann/DML/dml_links_author_H.html

[32] N. Nielsen, Die Gammafunktion. Chelsea Publishing Company, Bronx and New York, 1965.

[33] N.E. Nörlund, Vorlesungen über Differenzenrechnung. Chelsea, 1954. http://math-doc.ujf-grenoble.fr/cgi-bin/linum?aun=001355 http://dz-srv1.sub.uni-goettingen.de/cache/browse/AuthorMathematicaMonograph,WorkContainedN1.html

[34] M. Omarjee, Euler's constant with integrals. http://ohkawa.cc.it-hiroshima.ac.jp/AoPS.pdf/Euler's%20constant.pdf

[35] O. Schlömilch, Note sur quelques intégrales définies. J. Reine Angew. Math., 33, 316-324, 1846. http://www.digizeitschriften.de/

[36] O. Schlömilch, Sur l'intégrale définie $\int_0^\infty \frac{e^{-x\theta}}{\theta^2+a^2}d\theta$. J. Reine Angew. Math., 33, 325-328, 1846. http://www.digizeitschriften.de/

[37] J. Sondow, A faster product for $\pi$ and a new integral for $\log\frac{\pi}{2}$. Amer. Math. Monthly 112 (2005) 729-734. arXiv:math/0401406 [pdf]

[38] H.M. Srivastava and J. Choi, Series Associated with the Zeta and Related Functions. Kluwer Academic Publishers, Dordrecht, the Netherlands, 2001.

[39] E.O Tuck, On positivity of Fourier transforms. Bull. Austral. Math. Soc., 74(1), 133-138, 2006. http://www.austms.org.au/Gazette/2005/Sep05/Tuck.pdf

[40] E.T. Whittaker and G.N. Watson, A Course of Modern Analysis: An Introduction to the General Theory of Infinite Processes and of Analytic Functions; With an Account of the Principal Transcendental Functions. Fourth Ed., Cambridge University Press, Cambridge, London and New York, 1963.

Wessex House,
Devizes Road,
Upavon,
Pewsey,
Wiltshire SN9 6DL